\documentclass[12pt,a4paper,leqno]{article}
\usepackage{amsmath}
\usepackage{amssymb}
\usepackage{amsfonts}
 \usepackage{amsthm}

\usepackage[english]{babel}
\usepackage{mathrsfs}

\usepackage[applemac]{inputenc}
\usepackage{epigraph}
\numberwithin{equation}{section}
\renewcommand{\phi}{\varphi}
\renewcommand{\tilde}{\widetilde}
\usepackage[textwidth=15cm,  hoffset=-1.5cm, marginpar=4cm]{geometry} 

\DeclareMathOperator{\Ker}{Ker}

\DeclareMathOperator{\Id}{Id}

\def \R {\mathbb R}

\def \Ham {\mathscr H}
\def \DHam {\mathscr {DH}}

\begin{document}
\newtheorem{theorem}{Theorem}[section]
\newtheorem{lemma}[theorem]{Lemma}
\newtheorem{proposition}[theorem]{Proposition}
\newtheorem{Cor}[theorem]{Corollary}

\def \cal {\mathcal}

\theoremstyle{definition}
\newtheorem{definition}[theorem]{Definition}
\newtheorem{example}[theorem]{Example}
\newtheorem{exer}[section]{Exercise}
\newtheorem{conj}{Conjecture}
  \newtheorem*{def1.1}{Definition \ref{C0-comm}}

\theoremstyle{remark}
\newtheorem{remark}[theorem]{Remark}
\newtheorem{remarks}[theorem]{Remarks}
\newtheorem{question}{Question}
\epigraphwidth=4in
\epigraphrule=0pt

\font\smallcaps=cmcsc10
\font\nome=cmr8

\title{\textbf {Commuting Hamiltonians and Hamilton-Jacobi multi-time equations}}
\author{  \\{ \sc Franco Cardin${\ }^{1}$\qquad Claude Viterbo${\ }^{2}$}\\
\\
${\ }^{1}$ Dipartimento di Matematica Pura ed Applicata\\
 via Trieste 63 - 35121 Padova, Italia\\ {\tt   cardin@math.unipd.it}
\\
\\${\ }^{2}$  Centre de Math\'ematiques Laurent Schwartz \\ UMR 7640 du CNRS \\École Polytechnique - 91128 Palaiseau, France
\\ {\tt viterbo@math.polytechnique.fr   } }
\maketitle
\begin{abstract}The aim of this paper is twofold: first of all, we show that  the $C^0$ limit of a pair of commuting Hamiltonians commute. This means on one hand that if the limit of the Hamiltonians is smooth, the Poisson bracket of their limit still vanishes, and on the other hand that we may define ``commutation'' for $C^0$ functions. 

The second part of the paper deals with solving ``multi-time''  Hamilton-Jacobi equations using variational solutions. This extends the work of Barles and Tourin in the viscosity case to include the case of $C^0$ Hamiltonians, and removes their  convexity assumption, provided we are in the framework of ``variational solutions".\end{abstract}
\epigraph{Hamilton's variation principle can be shown to correspond to Fermat's {\it Principle} for a wave propagation in configuration space ($q$-space), and the Hamilton-Jacobi equation expresses Huygens' {\it Principle} for this wave propagation. Unfortunately this powerful and momentous conception of Hamilton, is deprived, in most modern reproductions, of its beautiful raiment as a superfluous accessory, in favour of a more colourless representation of the analytical correspondence. }{\textit {E. Schr\"odinger}, Quantization as a Problem of eigenvalues (Part II), Annalen der Physik, 1926}

\newpage

\section{Introduction}\label{intro}

The problem of finding solutions of multi-time Hamilton-Jacobi equations, by which one usually means  equations of the following type, where $x$ is in $\R^n$ and $t_j$ in $\R$ 
\begin{equation}\tag{MHJ}\left \lbrace
\begin{array}{ll}\frac{\partial}{\partial t_1}u(t_1,\ldots  ,t_d,x)+&H_1(t_1,\ldots ,t_d, x,\frac{
\partial}{\partial x}u(t_1, \ldots ,t_d,x))=0
\\ &\vdots \\ \frac{\partial}{\partial t_d}u(t_1,\ldots  ,t_d,x)+&H_d(t_1,\ldots ,t_d, x, \frac{ \partial
}{\partial x}u(t_1, \ldots ,t_d,x))=0
\end{array} \right . \end{equation}  with initial condition
$$u(0,...,0,x)=f(x)$$ has been initiated by Rochet in relation with some questions in economy, then by Lions-Rochet and studied more recently by Barles-Tourin
and Motta-Rampaz\-zo (\cite{Rochet,Lions-Rochet, Barles-Tourin,
Motta-Rampazzo}). 
Such a system of  equations is well-known to be overdetermined, and in order
to have a solution, we
need  the Hamiltonians to commute in a suitable sense. This is already obvious when applying  the method of
characteristics. Besides a suitable commutation condition, we need  to address the
question of the type of solution one is looking for. For first order equations, it is well
known, and was proved for more general equations
by Dacorogna and Marcellini \cite{Dacorogna-Marcellini}, that there
are plenty of $C^{0}$ solutions for such equations\footnote{By $C^0$ solution we mean $C^1$ almost everywhere, and satisfying the equation a.e.}. We then need
to select a particular ``class'' of solution, deemed to be the best suited to our problem.  A classical
choice is to look for viscosity solutions, which are the ``right'' solutions for optimal control and this is the type of
solution considered in the above papers (except for
\cite{Dacorogna-Marcellini}).

In Barles and Tourin's  paper, the existence of a viscosity solution is proved under the assumption that one of the
Hamiltonians is coercive with controlled growth\footnote{\label{fn1}Condition $(H1)$ in \cite{Barles-Tourin} states  that for each $R$ there exists $K_R$ such that

$ \vert H_1(x,p) \vert \leq K_R$ and 
$\vert \frac{\partial H_1}{\partial p}(x,p) \vert \leq K_R(1+ \vert x \vert) $
 in
 ${\mathbb R}^n\times B(0,R)$. 
Moreover the 
authors
assume $du_0(x)$ to be bounded. These conditions 
are stronger than the assumptions we need here.  
See appendix \ref{appendix-noncompact} for more
details. } 

(see \cite{Barles-Tourin}, page 1526, conditions (H1),(H2)), and more importantly 
that $H_{1},H_{2}$ satisfy the following conditions

  \begin{enumerate}
\item independent of $(t_1,...,t_d)$
\item  convex in $p$
\item  of class $C^1$ and satisfying the commuting condition
$\{H_j,H_k\}=0$
 \end{enumerate}

In fact the third condition
 can be weakened, as Barles and Tourin point out, to assume that the $H_j$ are $C^0$ and there are sequences of
Hamiltonians
$H_j^{\nu}$ of class  $C^1$   such that
  \begin{enumerate}
\item $\lim_{\nu\to \infty}H_j^{\nu}=H_j$ in the $C^0$ topology
\item \label{b} $\{H_j^{\nu},H_l^{\nu}\}=0$
\end{enumerate}

In other words, the $H_j$ are limits of commuting Hamiltonians.

However such an assumption is quite unpractical since it is already difficult to write two  commuting  $C^1$ Hamiltonians as nontrivial limits of commuting Hamiltonians.

The present paper has several goals.

First we solve the multi-time Hamilton-Jacobi equation in the framework of  ``variational solutions''
defined by Sikorav, Chaperon and the second author in \cite{Sikorav-pc, Chaperon, Viterbo-Ottolenghi} (see definition \ref{var-sol}). According to a result by Zhukovskaia (cf \cite{Zhukovskaia}),
if the  Hamiltonian is convex in $p$, the variational  solution must  coincide with the viscosity solution defined by Crandall and Lions (in \cite{Crandall-Lions}, see also our definition \ref{visc-sol-def}), so that our results extend those of Barles and Tourin. However in general these two solutions do not coincide (see an example in \cite{Viterbo-X}). For
variational solutions, we prove that we only need the $H_j$ to be
$C^0$, and condition \ref{b} is then replaced by a ``commutation''
condition best expressed in terms of symplectic invariants,
refining the following

\begin{definition}
\label{C0-comm} Let $H,K$ be two autonomous $C^0$ Hamiltonians. We shall say that
 $H$ and $K$ {\bf $\bf C^0$-commute}, if and only if there are sequences
 $H_\nu,K_\nu$ of $C^1$ Hamiltonians such that, all limits being for the $C^0$ topology, we have:
 \begin{enumerate}

\item $\lim_{\nu\to \infty} H_\nu=H$, $\lim_{\nu\to \infty} K_\nu=K$.

\item   $\lim_{\nu\to \infty} \{H_\nu,K_\nu\}=0$.
\end{enumerate} \end{definition}

A similar definition is given in section 3 for $H,K$ time-dependent and in
appendix \ref{Appendix-B} for the case of equations depending on
the function.

Postponing to the next Section the detailed geometrical setting for our Hamiltonians, the two main theorems of this paper are:

  \begin{theorem} \label{sec-Comm-Ham} If two $C^{1,1}$ Hamiltonians {\bf $\bf C^0$-commute}, then they commute in the usual
sense (i.e. their Poisson bracket vanishes).
 \end{theorem}
Here $C^{1,1}$ means differentiable with Lipschitz differential in the variables $(x,p)$. 
The above theorem tells us that our definition of $C^0$-commutation coincides with the classical one for smooth Hamiltonians. Note that this may be extended to the time-dependent setting as we shall see in section \ref{non-auton}. 

The above theorem sounds like a generalization of Eliashberg-Gromov's theorem 
(\cite{Gromov}, \cite{Eliashberg}, \cite{Ekeland-Hofer}) on the $C^0$ closure of the group 
of symplectic diffeomorphism, according to which the set of $2n$-tuples of functions 
on $\R^{2n}$, $(f_1,...,f_n,g_1,...,g_n)$ such that 
$$\{f_i,f_j\}=\{g_i,g_j\}=0,\; \{f_i,g_j\}=\delta_i^j$$ is 
closed\footnote{To be rigorous, Gromov and Eliashberg need  
the map $(x,y) \longrightarrow (f_1(x,y),..,f_n(x,y), g_1(x,y),...,g_n(x,y))$ 
to be bijective.}. We refer to \cite{Humiliere2} for an approach along these lines of the Gromov-Eliashberg theorem.  We also refer to improvemets of the above result from a quantitative point of view due to \cite{Entov-Polterovich-Zapolsky} using quasi-states. 

This has been extended by V. Humilière to other relations derived from so-called quasi-representations of finite dimensional Lie groups in the Poisson algebra. For example the Heisenberg relation $\{f,g\}=h; \{f,h\}=\{g,h\}=0$ is also $C^0$ closed (see \cite{Humiliere2}). 
 \begin{theorem} \label{main-thm} Assume the  Hamiltonians $H_1(t_1,..,t_d,x,p),...,H_d(t_1,..,t_d,x,p)$ on $T  {\mathbb R} ^n)$ satisfy the following conditions
 \begin{enumerate} 
 \item they are  locally Lipschitz in $(x,p)$ and their Lipschitz constant on the ball of radius $r$ has at most linear growth in $r$.
\item their support has an $x$-projection contained in a compact set.
\item They   {\bf $\bf C^0$-commute}. 
\end{enumerate} 
Then equation \thetag{MHJ} has a unique solution which is a variational solution of each individual equation.  If all  the Hamiltonians $H_j$'s are convex in
$p$, then $u$ is a viscosity  solution of each individual equation.
 \end{theorem}

  \begin{remarks}  \begin{enumerate} \item
  We refer to definition \ref{var-sol} and \ref{visc-sol-def} for the meaning of variational and viscosity solution. 
  
  \item  The growth condition is only needed to guarantee the existence of the flow for the
approximating Hamiltonians. Since if $dH$ has linear growth, we may
approximate it by smooth Hamiltonians $H^{\nu}$ such that the vector field  $X_{H^{\nu}}$
is complete\footnote{that is, the flow $\phi^t$ is defined for all $t$ in ${\mathbb R}$},
 this condition is sufficient. In fact, it is
enough to assume  there are constants $A,B$ such that  $ \vert
X_{H}(q,p) \vert  \leq A (\vert q \vert + \vert p \vert )+B$, i.e.
the Lipschitz norm of $H$ on a ball  grows at most linearly with
the radius. Note also that in some cases, we may guarantee
existence of the flow of the $X_{H^\nu}$ for other reasons. For
example when $H$ is  autonomous, and proper, since the
$H^{\nu}$ will satisfy the same assumption, and conservation of
energy implies that the flow remains in a compact set, the flow of $X_H$ is thus
defined for all times.
 \end{enumerate} \end{remarks}

 \begin{remark}
 The paper is not supposed to be written for specialists in symplectic topology, 
 although a certain familiarity with the basic constructions of \cite{Viterbo-STAGGF} 
 and \cite{Viterbo-Ottolenghi} is recommended. Appendices A and B are of a more symplectic 
 flavor and really address the question of Hamilton-Jacobi equations from a symplectic 
 topology viewpoint. In particular Appendix \ref{appendix-noncompact} addresses the question 
 of the growth conditions one must impose on the Hamiltonian and the  initial condition from  
 a purely geometric point of view, while Appendix \ref{Appendix-B} extends the main theorem 
 for equations depending on the unknown function.  \end{remark}

 \subsection{Organization of the paper} 
 Section 2 is devoted to a summary of the applications of Generating function theory to symplectic topology, in a slightly modified version with respect to \cite{Viterbo-STAGGF}. We also state the main properties of variational solutions as in \cite{Viterbo-Ottolenghi}.
 
 Section 3 deals with the proof of theorem 1.2. The proof is based on continuity  properties of the symplectic norm $c$ defined in Section 2. We prove that if $\{H_n,K_n\}$ are $C^0$, small, the flows $\phi_n^t, \psi_n^s$ of $H_n$ and $K_n$ have the following properties: 
 
 on one hand $t \longrightarrow \phi_n^t\psi_n^s \phi_n^{-t}\psi_n^{-s}$ is generated by a $C^0$ small Hamiltonian, and the properties of $c$ established in the previous section imply that  $c (\phi_n^t\psi_n^s \phi_n^{-t}\psi_n^{-s})$ is small. On the other hand, 
 if $H_n$ goes to $H$ with flow $ \phi^t$ and $K_n$ goes to $K$ with flow $\psi^t$, 
 $(\phi_n^t\psi_n^s \phi_n^{-t}\psi_n^{-s})$ goes to  $\phi^t\psi^s \phi^{-t}\psi^{-s}$. Uniqueness of limits and the fact that $c$ only vanishes on the identity implies that $\phi^t\psi^s \phi^{-t}\psi^{-s}=\Id$ for all $s,t$, hence $H$ and $K$ commute. 
 
In  section 4 we first show how multi-time equations have natural variational solutions, provided the Hamiltonians commute.  We then address a number of technical questions,  replacing the  invariants of section 2 by their stabilization. 

Section 5 eventually completes the proof of theorem 1.3. It is sufficient to deal with equation \thetag{MHJ} in the case of two Hamiltonians. Assume the two Hamiltonians are such that  $\{H_1,H_2\} $ is $C^0$  small. We then construct two Lagrangians, $L_{1,2}$ and $L_{2,1}$ obtained by ``solving''  the first of the two equations (for $t_2=0$) and then the second, and vice versa. We must then prove that the two Lagrangians $L_{1,2}$ and $L_{2,1}$ are close with respect to the $\gamma$ distance defined in \cite{Viterbo-STAGGF} and also in subsection 2.1(this is not so with respect to the $C^0$ distance). Once this is granted, it implies that the associated function $u_{1,2}$ and $u_{2,1}$ are $C^0$ close. 
The proof of the theorem is now obtained by limiting arguments. 

Appendix A gives the proof of some technical results. Appendix B extends the scope of the main theorems to the case of a non compact support. Appendix C,D,E give some complements on equations involving the unknown function, the geometric theory of Hamilton-Jacobi equations associated to coisotropic manifolds and historical comments. 

  \subsection{Acknowledgements}
  
 The authors warmly thank Franco Rampazzo for attracting their attention to this problem during a conference in Cortona, for
communicating his lecture notes on his work \cite{Motta-Rampazzo} and for many interesting discussions. We also would like to thank F.Camilli, I.Capuzzo Dolcetta and A.Siconolfi for the superb organization of the Cortona conference. Even though it is probably not related, we also mention the paper \cite{Rampazzo-Sussmann} on commutation of Lipschitz vector fields.

\section{Preliminary material}\label{prelim}

\subsection{Generating functions  and variational solutions of Hamilton-Jacobi equations}
 {\it We shall here assume that $N$ is a connected manifold without boundary, and either compact or that all Hamiltonians are
compact supported, as in \cite{Viterbo-STAGGF}. However, we shall explain in Appendix \ref{appendix-noncompact},
how our results extend to non-compact situations, provided we have some estimate on the growth of the Hamiltonians. }

Let $T^*N$ be the cotangent bundle of the manifold $N$ endowed with the canonical symplectic structure $\sigma=
\sum_{j=1}^n dp_j\wedge dx^{j}$.  To any Hamiltonian $H(t,z)$ (where $z=(x,p)$) on $ {\mathbb R} \times T^*N$ we associate the time-dependent  vector field $X_H$ defined by
$$\sigma(X_H,\xi)=-d_zH(t,z)\xi$$
 and the corresponding Hamiltonian flow $\phi^t_{t_{0}}$  defined by

\begin{gather*} \left\{ \begin{array}{ll} \frac{d}{dt}\phi^t_{t_{0}}=X_{H}(t,\phi^t_{t_{0}}) \\ \phi^{t_0}_{t_{0}}=Id \end{array} \right .
\end{gather*}

Let
$L$ be a Lagrangian submanifold of
$T^*N$ obtained   from the zero section
$0_N=\{(x,0)
\in T^*N
\mid x
\in N\}$  by the Hamiltonian isotopy $\phi^t_{t_{0}}$. We shall always assume that $\phi_{t_0}^t$ is well
 defined on $0_N$ for all $t$.  Then according to \cite{Sikorav} (relying on joint work  with Laudenbach in \cite{Laudenbach-Sikorav}),
 there exists a generating function quadratic at infinity  for $L$
i.e.  there exists a smooth function $S:N\times {\mathbb R} ^k \to
{\mathbb R} $ such that
\begin{enumerate}

\item $ (x,\xi) \mapsto \frac{ \partial S}{ \partial \xi} (x,\xi)$
has $0$ as a regular value (a Morse family in the terminology
introduced by  A. Weinstein in \cite{Weinstein-CBMS})
\item $S(x,\xi)=Q(\xi)$ for $\xi$ large enough, where $Q$ is a
non-degenerate quadratic form\footnote{This is sometimes
conveniently replaced by the condition $ \vert S(x,\cdot )-Q(\cdot
) \vert _{C^1}<+\infty$. We shall use these conditions
interchangeably in the rest of the paper. It is easy to prove that
existence of a  generating function of one kind is equivalent to
the  existence of a generating function of the other kind, see \cite{Brunella}. }.
\item
$$L=\{(x, \frac{ \partial S}{ \partial x} (x,\xi)) \mid \frac{ \partial S}{
\partial \xi}(x,\xi)=0\}$$
\end{enumerate}

In particular the critical points of $S$ are in one to one
correspondence with the points of $L\cap 0_N$.

  In the rest of the paper, we shall shorten the expression ``generating function quadratic at infinity'' by ``GFQI".

Let  $S^ \lambda =\{(x,\xi) \mid S(x,\xi) \leq \lambda \}$,
$E^\pm$ be the positive and negative eigenspaces of
$Q$,
 and $D^\pm$ be large discs in $E^\pm$. Since  for $c$ large enough $S^{\pm c}=N\times Q^{\pm c}$
we have

$$H^*(S^c,S^{-c})=H^*(N\times Q^c,N\times Q^{-c})=H^*(N)\otimes H^*(D^-, \partial D^-)$$
so that,   to each cohomology class $\alpha \in H^*(N)$ we may associate a class, image of $\alpha$ by the Künneth isomorphism, denoted by $T \alpha $. To  the class $T\alpha$ in
$H^*(S^c,S^{-c})$, we may associate a minimax critical level

$$c(\alpha,S)=\inf\{ \lambda \mid T \alpha \notin \Ker(H^*(S^c,S^{-c})
\to H^*(S^\lambda ,S^{-c})) \}
$$
Now it is proved in \cite{ Theret} and \cite{Viterbo-STAGGF} 
 that given $L$, $S$ is essentially unique up to adding a constant, and more precisely,  up to a global shift, the numbers $c(\alpha , S)$ depend only on $L$, not on $S$, and they are thus
denoted by $c(\alpha,L)$. Moreover, denoting by $$\gamma (L)=c(\mu,L)-c(1,L)$$ where $1\in H^0(N), \mu \in
H^n(N)$ are generators.  We know that $\gamma(L)$ is well defined and vanishes  if and only if
$L=0_N$ (see for instance \cite{Viterbo-STAGGF}) .

Moreover let $S_{x}(\xi)=S(x,\xi)$ be the restriction of $S$  to the fiber over $x$. We can look for a minimax
as above for the function  $S_x$. Since the cohomology of the point is one dimensional, denoting
its  generator by $1_{x}$,  we set

\begin{definition}\label{var-sol} The continuous function
$$u_L(x)=c(1_{x},S_{x})$$ is called {\bf the variational solution} of the equation $$(x,du(x))\in L.$$ In particular if $L\subset H^{-1}(0)$ $u_L$ is  
 {\bf a variational solution} of the Hamilton-Jacobi equation $$H(x,du(x))=0$$ 
 \end{definition} 
It has been proved by  Sikorav and Chaperon (see \cite{Chaperon, Sikorav-pc} and also \cite{Viterbo-Ottolenghi}), that such a function is indeed a solution
of the {\it Lagrangian Hamilton-Jacobi equation}, that is
$$
(x,du_L(x))\;\text{is in }\; L \;\text{for almost all $x$ in $N$}
$$
When  $L\subset H^{-1}(0)$, we have a solution of the classical Hamilton-Jacobi equation 

 $$H(x,du_L(x))=0 \;\text{for almost all $x$  in $N$}$$
Of course, for any constant $c$, $u_L(x)+c$ is also a solution.
For evolution equations, that is

 \begin{gather*} \left \{ \begin{array}{ll} \frac{\partial u}{\partial t}(t,x)+H(t,x,\frac{
\partial u}{\partial x}(t,x))=0 \\ u(0,x)=f(x) \end{array}\right .   \end{gather*}
the construction of variational solutions for a single equation can be rephrased  as follows.
Let $\Lambda_0$ be a Lagrangian submanifold of $T^*N$, and $H(t,z)$ be  
a time-dependent Hamiltonian. We consider $\tilde \Lambda_{0,H}=\{(0,-H(0,z),z) \mid z \in \Lambda_0\}$, and
$$L= \bigcup_{t\in \R} \Phi^t(\Lambda_{0,H})\subset T^*(\R \times N)$$

In \cite{Viterbo-Ottolenghi}, it is also proved that for evolution
equations  the definition of variational solutions extends to $C^0$
Hamiltonians:  indeed if $H^{\nu}$ tends to $H$, then the solution $u^\nu$
converges to $u$. This follows from the property $$ \vert u_H-u_K \vert_{C^0([0,T]\times \R^{n})} \leq T \Vert H-K \Vert _{C^0([0,T]\times \R^{2n})}$$ which in turn follows from Proposition \ref{c-conv} (see \cite{Viterbo-Ottolenghi, Viterbo-X}). Note that in the framework of viscosity
solutions, this property is called {\bf stability}. A priori, even though $u$ is Lipschitz -hence according to
Rademacher theorem is almost everywhere differentiable- we do not
claim that $u$ satisfies the Hamilton-Jacobi equation almost everywhere\footnote{This is the case for viscosity solutions, but is unknown for variational solutions.}.

\subsection{Capacities for Hamiltonian flows}

Let  $L_1,L_2$ be Lagrangian submanifolds generated by
$S_1,S_2$.  We may define  $c(\alpha, S_1-S_2)$, and as this does not depend on the choice
 of $S_1,S_2$ but only on $L_1,L_2$ we denote it by abuse of language $c( \alpha ,L_1-L_2)$, even though it is
not really determined by the set $$L_1-L_2=\{(x,p_1-p_2) \mid (x,p_1) \in L_1, (x,p_2)\in L_2\}$$ but depends on both
$L_1$ and $L_2$.

 We denote by $\gamma (L_1-L_2)$ the difference $c(\mu, L_1-L_2)-c(1,L_1-L_2)$. This is non-negative according to \cite{Viterbo-STAGGF} and vanishes if and only if $L_1=L_2$.

 Now to a compact supported symplectic isotopy, denoted by $\psi$,
we may associate a symplectic invariant as follows:

 \begin{definition}

Let ${\mathscr L}$ be the set of Lagrangian submanifolds Hamiltonianly isotopic to the zero section,    $\Ham (T^*N)$ be the set of smooth time-dependent Hamiltonians on $ {\mathbb R} \times T^*N$, $\DHam (T^*N)$ be  the set of time one maps of such Hamiltonians. We shall use the notation $\Ham, \DHam$ if there is no ambiguity.
We set for $\psi \in \DHam $,
$$\tilde\gamma(\psi)=\sup\{\gamma(\psi(L)-L) \mid L\in {{\mathscr L}} \}$$
 \end{definition}

 We shall need the following

 \begin{proposition}\label{2.3}
 \begin{enumerate}
\item $$\tilde\gamma(\psi)\geq 0$$ and $\tilde\gamma(\psi)=0$ if
and only if $\psi=\Id $ \item $$\tilde\gamma
(\psi^{-1})=\tilde\gamma(\psi)$$ \item (triangle inequality)
$$\tilde\gamma(\psi\varphi)\leq \tilde\gamma(\psi)+\tilde\gamma
(\varphi)$$ \item (invariance by conjugation)
$$\tilde\gamma(\varphi \psi \varphi^{-1})=\tilde\gamma(\psi)$$
\end{enumerate}
 \end{proposition}
 \begin{proof}

The proof of this proposition is postponed to Appendix \ref{Pf}
\end{proof}

 \begin{definition} \label{def-c-conv} We shall say that the sequence $\phi_n$ in $\DHam$  {\bf
$c$-converges} to $\phi$ if and only if $$\lim_{n\to
\infty}\tilde\gamma(\phi_n^{-1}\phi)=0$$
We shall use the notation $$\phi_n\overset{c}\longrightarrow \phi$$ for $c$-convergence\footnote{$c$ stands for capacity, see \cite{Ekeland-Hofer}.}.
 \end{definition}
  \begin{remark}
  Since our invariant is called $\gamma$, we should talk about $\gamma$-convergence. In fact our $c$-convergence is indeed related to $\Gamma$-convergence in the calculus of variations, but we want to avoid any confusion here.
  \end{remark}
 We also need the following estimate:

  \begin{proposition}\label{c-conv} Assume $\psi$ is the time one map of the Hamiltonian $H(t,z)$. Then  we have
$$\tilde\gamma(\psi) \leq \Vert  H \Vert_{C^0}:=\sup_{(t,z)\in [0,1]\times {\mathbb R}^{2n}} H(t,z) -\inf_{(t,z) \in [0,1]\times {\mathbb R}^{2n}} H(t,z)$$  As a consequence, if $\phi_n$ and $\phi$ are generated by
$H_n$ and $H$, and  if $H_n\longrightarrow H$ in the $C^0$ topology, then $\phi_n$
$c$-converges to $\phi$.
 \end{proposition}
 \begin{proof}
We first prove that if $H_0, H_1$  are Hamiltonians and
$\phi_0,\phi_1$ are their flows,  we may normalize $S_0, S_1$
generating functions of $\phi_0(L),\phi_1(L)$ so that their critical
values are  those of
$$
{\mathcal S}_{j}(x(1);  \gamma,\xi)=S(x(0),\xi)+\int_0^1 p(t)\dot x(t) +H_j(t,x(t),p(t))dt
$$
with respect to the (infinite) auxiliary parameters $\gamma,\xi$,  and were $S$ is generating
function quadratic at infinity for $L$, and $\gamma=(x,p):[0,1]\to T^{*}N$.
 Indeed we have

 \begin{gather*}D{\mathcal S}_j(x(1); \gamma,\xi) \left(\delta x(1), \delta x(0),\delta \gamma, \delta \xi\right) = 
 \\ 
 p(1)\delta x(1) +
 \frac {\partial S}{\partial \xi}(x(0),\xi) \delta\xi - \left ( p(0)-\frac{\partial S}{\partial x}(x(0),\xi)\right )\delta x(0) +\\  \int_0^1 \left [ \left ( \dot x (t)+ \frac{\partial H_j}{\partial p}(x(t),p(t))\right ) \delta p(t) - \left ( \dot p (t)-\frac{\partial H_j}{\partial x}(x(t),p(t))\right ) \delta x(t) \right ] dt \end{gather*}

 According to \cite{Viterbo-STAGGF}, \cite{Theret} generating functions associated to a Lagrangian are
``essentially unique'', up to a constant. Thus the critical values
of two such functions differ by a global translation. Since
${\mathcal S}_{j}(x(1);\gamma,\xi)$ is formally a generating function,
and in particular has critical values coinciding up to translation
with those of any other generating function, we may use its
critical values to normalize the $S_{j}$ (i.e. we replace $S_{j}$
by $S_{j}+c_{j}$ so that the critical values of $S_{j}+c_{j}$ and
$\mathcal { S}_{j}$ coincide).

  In particular we claim that if $H_0\leq H_1$  we have $c(\alpha, S_0)\leq
  c(\alpha,S_1)$.

For this we argue as in the proof of proposition 4.6 from \cite{Viterbo-STAGGF}.  We consider the linear
interpolation
$H_\lambda(t,x,p)=(1-\lambda) H_0(t,x,p)+\lambda H_1(t,x,p)$. Let
$\phi_{\lambda}^t$ be the flow of $H_\lambda$ and $\phi_\lambda$ be the time one map. The associated generating function
$S_\lambda (x,\xi)$ of
$\phi_{\lambda}(L)$ is normalized as above. Now for a critical point of
$S_\lambda$, we get an intersection point $(x_\lambda,0)=(x_\lambda(1),p_\lambda(1))$ of $\phi_{\lambda}(L)\cap 0_N$ and the critical values are  those of
$$
{\mathcal S}_\lambda(x(1);  \gamma,\xi)= S(x(0),\xi)+\int_0^1 [p(t)\dot x(t)+H_\lambda(t,x(t),p(t))]dt
$$
corresponding to critical points of the form $(x_\lambda(1),\gamma_\lambda, \xi_\lambda)$ where $$\gamma_\lambda(t)=(x_\lambda(t),p_\lambda(t))=\phi_{\lambda}^t(x_\lambda(0),p_\lambda(0))$$
and $$\frac{\partial}{\partial \xi}S(x_\lambda(0),\xi_\lambda)=0$$

Using the fact that generically  $x_\lambda=x_\lambda(1)$ is piecewise $C^1$, and $S(x_\lambda,\xi_\lambda)$ is continuous, it is enough to know that for all but a finite set of values of $\lambda$ we may write
$$\frac{d}{d\lambda}S_\lambda(x_\lambda,\xi_\lambda)= \int_0^1\frac{d}{d\lambda}H_\lambda(t,x_\lambda(t),p_\lambda(t))$$

and this quantity is positive if $H_0\leq H_1$.

Now  $ \Vert H \Vert_{C^0} = C$ means that $a\leq H \leq b$ with $b-a\leq C$. Then, since for the constant
Hamiltonian $h_a(x)=a$, $$c(\mu,S_a)=c(1,S_a)=a$$ (again because of the above normalization) and we get $a\leq
c(1,S_H)\leq c(\mu,S_H)\leq  b$ and  we have $\gamma(S_H)\leq b-a= \Vert H \Vert _{C^0}$ 
 \end{proof}
\begin{remark} As the referee pointed out, A. Weinstein noticed long ago that the action functional is a "generating function" in some generalized sense. Taking some finite dimensional reduction of this, one can associate to $H$ a GFQI $S_H$ such that  $$\Vert S_H-S_K \Vert \leq \Vert H-K \Vert.$$
\end{remark} 

\begin{Cor}\label{Cor-c-conv}
Assume $H_\nu \longrightarrow H$ in the $C^0$ topology, where $H$ is
in $C^{1,1}$, (i.e. it has  Lipschitz derivatives), and $\phi_\nu^t, \phi^t$ are the flows of $H_\nu$ and
$H$. Then for all $t$,
$\phi_\nu^t$ c-converges to $\phi^t$:
$\phi_\nu^t \overset{c}\longrightarrow
\phi^t$.
\end{Cor}

\begin{proof}
Indeed  $\phi^{-t}\phi_\nu^t$ is the flow of
$H_\nu(t,\phi^t(z))-H(t,\phi^t(z))$, and clearly if $H_\nu
\longrightarrow H$ in the $C^0$ topology, this quantity goes to
zero, hence $\tilde\gamma (\phi^{-t}\phi_\nu^t)$  
goes to zero, which is
equivalent to $\phi_\nu^t \overset{c}\longrightarrow \phi^t$.
\end{proof}

\begin{remark}
Given $\phi$ the time one map of a symplectomorphism of the symplectic manifold $M$,  we may define its {\bf Hofer norm} as follows 
(\cite{Hofer})
Let ${\mathcal H}(\phi)$ be the set of (time-dependent) Hamiltonians on   $M$ such that
the time one flow associated to $H$ is $\phi$.

Then $$ \Vert \phi \Vert = \inf \left \{ \int_0^1 [ \max_x H(t,x)- \inf_x H(t,x) ] dt \mid H \in {\mathcal H}(\phi)\right \}$$

In fact, the proofs of both proposition 2.6 and corollary 2.7 show more than stated: we prove that the identity map from $({\mathcal H}, \Vert \bullet \Vert ) $ to $({\mathcal H}, \tilde\gamma )$ is a contraction, or in less pedantic terms, that $$ \tilde \gamma (\phi) \leq \Vert \phi \Vert $$
\end{remark}

This allows to set the following definition, as in \cite{Humiliere}:


\begin{definition} (\cite{Humiliere})
We define ${\mathfrak H}_\gamma (M)$ the completion of $\Ham (M)$ for the metric $\gamma$, and $\mathfrak{HD}_\gamma (M)$ the completion of $\DHam (M)$ for the metric $\tilde\gamma$.
\end{definition} 

Let us point out that the ``time one flow'' map $\Ham (M)\longrightarrow \DHam (M)$
extends to a continuous map $\mathfrak {H}_\gamma (M) \longrightarrow
\mathfrak{HD}_\gamma (M)$.
According to Proposition 2.6, we have a continuous map of $C^0( {\mathbb R} \times M) $
into  ${\mathfrak H}_\gamma (M)$. Moreover, according to a theorem of V.
Humili\`ere (\cite{Humiliere}), Hamiltonians with some controlled singularities
also live in ${\mathfrak H}_\gamma (M)$.

\subsection{Viscosity solutions}

The only fact the reader needs to know here about viscosity solutions is  Zhukovskaia's theorem.  For his convenience, we repeat  the definition of viscosity solutions in the framework of evolution equations 

\begin{definition} \label{visc-sol-def}
A viscosity subsolution (resp. supersolution) of the Hamilton-Jacobi equation
 \begin{gather*} \left \{ \begin{array}{ll} \frac{\partial u}{\partial t}(t,x)+H(t,x,\frac{
\partial u}{\partial x}(t,x))=0 \\ u(0,x)=f(x) \end{array}\right .   \end{gather*}
is a function $u$ satisfying the initial condition and such that if $\phi(t,x)$ is a function such that  $u(t,x)-\phi(t,x)$ has a local maximum (resp. minimum) at $(t_0,x_0)$ we have 

\begin{gather*}\begin{array}{ll} \frac{\partial u}{\partial t}(t,x)+H(t,x,\frac{
\partial u}{\partial x}(t,x))\leq 0  \;\; (\text{resp.}\; \geq 0) \end{array}  \end{gather*}
Moreover $u$ is a viscosity solution if and only if it is both a viscosity subsolution and a viscosity supersolution.
 \end{definition}
We now have 

 \begin{theorem} [Zhukovskaia's theorem] (\cite{Zhukovskaia}, \cite{Bernardi-Cardin})
 If $H$ is convex in $p$ then a function $u(t,x)$ is a viscosity solution of the Hamilton-Jacobi equation
  \begin{gather*} \left \{ \begin{array}{ll} \frac{\partial u}{\partial t}(t,x)+H(t,x,\frac{
\partial u}{\partial x}(t,x))=0 \\ u(0,x)=f(x) \end{array}\right .   \end{gather*}
if and only if it is a variational solution.
 \end{theorem} 

Note that an example was given in \cite{Viterbo-Ottolenghi}  showing that this theorem fails if we remove the convexity assumption. 

  \section{Commuting autonomous Hamiltonians and the proof of theorem \ref{sec-Comm-Ham}}

 Let  $f:P_{1}\to P_{2}$ be a diffeomorphism between two manfiolds,  and $X$ be a vector field on $M_{1}$, then
 the  {\it push-forward} of $X$,  $f_*X$ is defined as the vector
field: $f_*X(y)=df(f^{-1}(y))X(f^{-1}(y))$.
 Its main property  is that if  $\phi^t$ is 
 the flow of  $X$, then
$$\frac{d}{dt} (f \circ \phi^t)(x)=(f_*X)(f \circ \phi^t (x))$$
since $$\frac{d}{dt} (f\circ
\phi^t)(x)=df(\phi^t(x))X(\phi^t(x))=df(f^{-1} f\circ \phi^t (x)
)X(f^{-1} f
 \circ\phi^t(x))=f_* X(f \phi^t(x))$$

 \begin{proposition} \label{lemma-comm} Let $H,K$ be two autonomous $C^{1,1}$ Hamiltonians of the symplectic manfiold $M$ and assume
$\{H,K\}$ to be
$C^0$ small. Then denoting by $\phi^t, \psi^s$ the Hamiltonian flows of $H$ and  $K$,
the Hamiltonian isotopy $t \mapsto \phi^t\psi^s\phi^{-t}\psi^{-s}$ is generated by a $C^0$ small (time-dependent) Hamiltonian.
\end{proposition}

 \begin{proof}
 Indeed we have, setting
$u=\phi^t\psi^s\phi^{-t}\psi^{-s}(x)$:

\begin{gather*}
\frac{d}{dt}(\phi^t\psi^s\phi^{-t}\psi^{-s})(x) =
\frac{d}{dt}\phi^t(\psi^s\phi^{-t}\psi^{-s}(x)) +
d(\phi^t\psi^s)(\phi^{-t}\psi^{-s}(x))(\frac{d}{dt}\phi^{-t})(\psi^{-s}(x)) =\\
X_H(\phi^t\psi^s\phi^{-t}\psi^{-s}(x)) -  d(\phi^t\psi^s)(\phi^{-t}\psi^{-s}(x))X_H(\phi^{-t}\psi^{-s}(x))=\\
X_H(u)-d(\phi^t\psi^s)(\psi^{-s}\phi^{-t}(u))X_H(\psi^{-s}\phi^{-t}(u))=\\  [X_H-(\phi^t\psi^s)_*X_H](u),
\end{gather*}  and since for each symplectic diffeomorphism, $\rho$, we have
$\rho_*X_L=X_{L\rho^{-1}}$, the vector field
$[X_H-(\phi^{t}\psi^s)_*X_H]$ is Hamiltonian, with Hamiltonian function
$$ L_{s}(t,x)=H(x)-H (\psi^{-s} \phi^{-t}(x))
$$

We may thus compute
$$
\frac{\partial}{\partial s}L_{s}(t,x)=-\frac{d}{d s}H (\psi^{-s}\phi^{-t} (x))=dH(\psi^{-s}\phi^{-t}(x))\cdot
X_K(\psi^{-s}\phi^{-t}(x))=\{H,K\}(\psi^{-s}\phi^{-t}(x))
$$ and since  $H(\phi^{-t}(x))=H(x)$,  we have $L_{0}=0$ and we may then estimate
 \begin{gather*}
\vert L_{s}(t,x) \vert = \vert  L_{s}(t,x)-L_{0}(t,x)\vert \leq \int_0^s \vert
\frac{\partial}{\partial \sigma}L_{\sigma}(t,x)\vert d\sigma \leq \\ \int_0^s\vert
\{H,K\}(\psi^{-\sigma} \phi^{-t}(x)) \vert d\sigma \leq s\Vert \{H,K\}\Vert_{C^0}
  \end{gather*}

Thus if $\vert
\{H,K\} \vert$ is $C^0$  small, for each $s$,  the flow $t\mapsto
\phi^t\psi^s\phi^{-t}\psi^{-s}$ is generated by a $C^0$  small Hamiltonian. \end{proof}

Remember that according to definition \ref{C0-comm}, two autonomous Hamiltonians  $H$ and $K$ {\bf
 $\bf C^0$-commute} if and only if there exist
  sequences of smooth Hamiltonians
$H_n,K_n$ such that, in the $C^0$ topology:  $H_n$ goes to $H$, $K_n$ goes to $K$, 
and  $\{H_n,K_n\}$ goes to zero.

Fortunately this definition does not conflict with the standard one according to theorem \ref{sec-Comm-Ham} that we now prove:
\begin{proof}[Proof of theorem \ref{sec-Comm-Ham}]  Remember  that all Hamiltonians are compact supported. We shall explain in appendix \ref{appendix-noncompact}
how to extend this to more general situations.

We also assume temporarily that $H,K$ are  of class $C^{1,1}$.
 Let now $H_n,K_n$ be a sequence of compact
supported $C^{1,1}$ Hamiltonians, such that their $C^0$-limits satisfy:

\begin{enumerate}
\item $\lim_{n\to \infty} H_n = H, \lim_{n\to \infty} K_n = K$

\item
$\lim_{n\to \infty} \{H_n,K_n\}=0$
\end{enumerate}

We wish to prove that if $H,K$ are $C^{1,1}$ then $\{H,K\}=0$. Let $\phi^t_n, \psi^s_n$ the flows of $H_n, K_n$,
and set
$\rho_n(s,t,x)=\phi_n^t\psi_n^s\phi_n^{-t}\psi_n^{-s}(x)$. The flow
$t\mapsto \rho_n(s,t,x)$ is the flow of the Hamiltonian $L^n_{s}(t,x)$.

Making use of the topology of $c$-convergence from  definition
\ref{def-c-conv} above and applying Corollary \ref{Cor-c-conv} we
have that $\phi^t_n$ $c$-converges to $\phi^t$, and $\psi^s_n$ to
$\psi^s$. Using the triangle inequality for $\tilde \gamma$ we
see that $\phi_n^t\psi_n^s\phi_n^{-t}\psi_n^{-s}$ $c$-converges to
$\phi^t\psi^s\phi^{-t}\psi^{-s}$. On the other hand, according to
proposition \ref{lemma-comm}, since $\{H_n,K_n\}$ is $C^0$ small,
$\tilde\gamma(\phi_n^t\psi_n^s\phi_n^{-t}\psi_n^{-s})$ goes to
zero and thus, using proposition \ref{c-conv},  we get that
$\phi^t\psi^s\phi^{-t}\psi^{-s}=\Id$. Since this holds for any $s,t$, it obviously implies that  $H$ and $K$ commute.
This concludes our proof.
 \end{proof}

The proof required the Hamiltonians to be $C^{1,1}$ in order to define the flows of $X_H,X_K$. One should be able to deal with the slightly more general $C^{1}$ case by using the methods of \cite{Humiliere}.

 The same proof shows that the following definition of commutation is also compatible with the standard one

 \begin{definition} \label{c-commute}
 Let  $H,K$ be two $C^0$ Hamiltonians on $T^*N$.  We shall say that $H$ and $K$ {\bf $c$-commute} if and only if there
exists sequences $H_n,K_n$ such that denoting the flows  of $H_n$ and $H$
 by $\phi^t_n, \phi^t$, and those of $K_n$  and $K$ by
$\psi_n^t, \psi^t$ (since $H,K$ are only $C^0$, the flows are only defined in $\mathfrak{HD}_\gamma (M)$), we have

$$\phi_n^t {\overset c \longrightarrow } \phi^t$$ $$\psi_n^t  {\overset c \longrightarrow }
\psi^t$$

$$\phi_n^t\psi_n^s\phi_n^{-t}\psi_n^{-s} {\overset c \longrightarrow}  Id$$

\end{definition}

 Note that even though this is a very natural notion, it is not the one we will need to generalize  theorem \ref{sec-Comm-Ham}.

\subsection{The time-dependent and multi-time case}\label{non-auton}
\subsubsection{The time-dependent case} Assume now $H, K$ are time-dependent, that is of the form
$H(t,z)$ and
$K(t,z)$. Commutation is then expressed by different formulas.  If $X_H,X_K$  are the Hamiltonian vector fields,
we define
$$\lbrack H,K\rbrack (t,z)= d_zH(t,z)X_K(t,z)+ \frac{\partial H}{\partial t}(t,z) -\frac{\partial K}{\partial
t}(t,z)
$$

In the sequel, even for time-dependent vector fields we denote by $\{H,K\}(t,z)$ the  Poisson bracket of $H,K$
for fixed time, that is $\{H,K\}(t,z)=d_xH(t,z)X_K(t,z)$, so that

$$\lbrack H,K\rbrack (t,z)= \{H,K\}(t,z)+ \frac{\partial H}{\partial t}(t,z) -\frac{\partial K}{\partial t}(t,z)
$$

The vanishing of
$\lbrack H,K\rbrack$ has the following interpretation. To $H(t,z)$ we may associate the autonomous Hamiltonian
$$\widetilde  H(t,\tau,z)=\tau+H(t,z)$$ defined on $({\mathbb R}^2 \times M ,  d\tau\wedge dt+ \omega)$, and
the  associated flow is
 $$\Phi^s:  (t,\tau, z) \mapsto (t+s, \tau+H(t,z)-H(t+s,\phi_t^{t+s}(z)),
\phi_t^{s+t}(z))$$

where $\phi_t^{t+s}$ is the flow on $M$  of the non autonomous Hamiltonian $H(t,x)$.

Similarly to $K$ we associate $\Psi^s$.  Now it is easy to check that  $[H,K](t,z)=\{\widetilde H, \widetilde  K\}(t,\tau,z)$, the
commutation of $\Phi^t$ and $\Psi^s$ implies  obviously that $$\psi_b^d\phi_a^{b}(z) = \phi_c^{d}\psi_a^c(z)$$
for any $a,b,c,d$ with $d-b=c-a$.

With this definition, the results of the previous section can be extended to the time-dependent case.

\begin{definition}\label{C0-comm-non-aut} Let $H(t,z),K(t,z)$  be two continuous Hamiltonians. We shall say that
$H$ and
$K$ $C^0${\bf -commute}
 if and only  if there are sequences
$H_n,K_n$ such that
$\lbrack H_n,K_n \rbrack$ goes to zero in  the  $C^0$ topology. 
\end{definition}

If $H,K$ are of class $C^{1,1}$  this definition
coincides with the condition $\lbrack H, K\rbrack =0$.

\subsubsection{The multi-time case} We shall limit ourselves to the case $d=2$. We consider $H(s,t,z), K(s,t,z)$ and we consider the
Hamiltonian flows
$$\frac{d}{dt} \phi^{t}_{t_{0}}(s;x)=X_{H}(s,t,\phi^{t}_{t_{0}}(s;x))$$

$$\frac{d}{ds} \psi^{s}_{s_{0}}(t;x)=X_{K}(s,t,\phi^{s}_{s_{0}}(t;x))$$

Then, we consider the Hamiltonians on $T^{*}(N\times {\mathbb R}^{2})$ given by
$$\widetilde H(s,t,\sigma,\tau,x)=\tau +H(s,t,z)$$
$$\widetilde K(s,t,\sigma,\tau,x)=\sigma +K(s,t,z)$$

Denote their flows by $\widetilde\Phi^a, \widetilde\Psi^b$ respectively. The last component of
$\widetilde\Phi^{a}(s,t,\sigma,\tau)$ is given by  $ \phi^{t+a}_{t}(s;z)$ and the last component of
$\widetilde\Psi^{b}(s,t,\sigma,\tau)$ by $\psi^{s+b}_{s}(t;z)$.

 Now $$\{\widetilde H , \widetilde K\}= \{H,K\}(s,t,z)+ \frac{\partial H}{\partial t}(s,t,z) -\frac{\partial
K}{\partial s}(s,t,z)$$

 and the commutation of the flows $\phi^{s}_{s_{0}}(t;x)$ and $\psi^{t}_{t_{0}}(s;x)$ is given by the vanishing
of the above expression.

 \begin{definition}
 We set $$\ll H,K\gg=\{H,K\}(s,t,z)+ \frac{\partial H}{\partial t}(s,t,z) -\frac{\partial K}{\partial s}(s,t,z)$$
  \end{definition}
   \begin{proposition}
   The ``flows'' $\phi^{s}_{s_{0}}(t;z)$ and $\psi^{t}_{t_{0}}(s;z)$ commute, or more explicitly
  $$\phi^{s}_{s_{0}}(t;\psi^{t}_{t_{0}}(s_{0};z))= \psi^{t}_{t_{0}}(s;\phi^{s}_{s_{0}}(t_{0};z))$$
   if and only if
   $\ll H,K\gg =0$.
    \end{proposition}

    \begin{remark}
    By applying Proposition \ref{lemma-comm}, we get that all our different notions of commutation (in the time-dependent or multi-time case) are well defined in the continuous case, and the definition is compatible with the usual definition in the $C^{1,1}$ case.
    \end{remark}
\section{Multi-time  Hamilton-Jacobi equations and the proof of theorem \ref{main-thm} }

Consider the  system of equations \thetag{MHJ} defined in the introduction. 
 For small $t_1,...,t_d$, and a Hamiltonian of class  $C^{1,1}$, we can find a solution to the above equation provided the Hamiltonians
commute in the sense of the previous section, that is the functions $\tau_j+H_j(t_1,..,t_d,x,p)$ commute on
$T^*({\mathbb R} ^d\times N )$, which we can rewrite as:

\begin{equation}\label{H1}  \{H_j,H_k\}(t,z)+ \frac{\partial H_j}{\partial t_k}(t,z) -\frac{\partial
H_k}{\partial t_j}(t,z)=0 \;\; t\in {\mathbb R}^d, z \in T^*N \tag{H1}
 \end{equation}

    \subsection{The smooth case} Under assumption (\ref{H1}) it is easy to check that the method of
characteristics will yield a solution for $t_j$ small.

Let us remind the reader of  the geometric approach to this method.
    Set
$$\Lambda_{f;H_1,...,H_d}=\{(0,...,0,x,-H_1(0,...,0,x, df(x)),... ,-H_d(0,...,x,df(x)),df(x))\mid x
\in N\}$$

This is an isotropic submanifold in $T^*(N\times {\mathbb R} ^d)$,
and  the vector fields $X_j$ associated to $\widetilde
H_j(t_1,...,t_d,\tau_1,...,\tau_d,
x,p)=\tau_j+H_j(t_1,..,t_d,x,p)$ being pairwise in involution, are all
 tangents to the submanifold $$\Sigma=\bigcap_{j=1}^d \widetilde H_j^{-1}(0)$$ and transverse to $\Lambda_{f;H_1,...,H_d}$.  As a result the
flows $\widetilde\Phi_j^t$ of $X_j$ commute and defining

$$L_{f;H_1,...,H_d}=\bigcup_{t_1,...,t_d \in {\mathbb R} ^d} \widetilde\Phi_1^{t_1} \ldots
\widetilde\Phi_j^{t_j} \ldots
\widetilde\Phi_d^{t_d} (\Lambda_{f;H_1,...,H_d})$$

the commutation property insures that the order in which  we take the $\widetilde \Phi_j^{t_j}$ is irrelevant, and
that  $L_{f;H_1,...,H_d}$ is a Lagrangian submanifold contained in $H_{j}^{-1}(0)$, hence in $\Sigma$ .

    The following lemma is a straightforward extension of a result from \cite{Viterbo-Ottolenghi}. 
 \begin{lemma} The manifold $L_{f;H_1,... , H_d}$ is Hamiltonianly isotopic
 to  the zero section in $T^*( {\mathbb R} ^d \times N)$.
\end{lemma}
    \begin{proof}

Consider the family $\Phi^{\tau}_{j,\lambda}$ of flows of the Hamiltonian $\tau_j+\lambda H(t_1,...,t_d,x,p)$.
We have $\Phi^{\tau}_{j,1}=\Phi^{\tau}_j$ and $\Phi^{\tau}_{j,0}$ is the translation by $\tau$ in the
$t_j$ coordinate.

Now the family  $\Phi^{\tau}_{j,\lambda}$  is a continuous family of proper maps, and the  family of Lagrangian
submanifolds given by

$$L^{\lambda}_{f;H_1,...,H_d}=\bigcup_{t_1,...,t_d \in {\mathbb R} ^d} \widetilde\Phi_{1,\lambda}^{t_1}
\ldots
\widetilde\Phi_{j,\lambda}^{t_j} \ldots
\widetilde\Phi_{d,\lambda}^{t_d} (\Lambda_{\lambda f;H_1,...,H_d})$$

is an exact proper Hamiltonian  isotopy of Lagrangian submanifolds.

Note that we do not need here that the Hamiltonians commute, and we shall indeed use the lemma in the more
general setting.
 \end{proof}

Now, still  assuming the $H_j$ are of class $C^{1,1}$ in $(x,p)$, section \ref{prelim} implies that $L_{f;H_1,... , H_d}$ has a GFQI. We thus have a variational solution
$u_{f,H_1,...,H_d}$ of the following equations

\begin{equation*}\tag{$\star$} \left \lbrace
\begin{array}{llll}\frac{\partial}{\partial t_1}u(t_1,t_2,\ldots  ,t_d,x)&+&H_1(t_1,t_2,\ldots ,t_d, x,
\frac{\partial}{\partial  x}u(t_1,t_2, \ldots ,t_{d-1},t_d,x))=0 \\
  \\
\frac{\partial}{\partial t_2}u(0,t_2,\ldots  ,t_d,x)&+&H_2(0,t_2\ldots ,t_d, x, \frac{ \partial }{\partial
x}u(0,t_2 \ldots ,t_{d-1},t_d,x))=0 \\
 \vdots &&\vdots \\
\frac{\partial}{\partial t_d}u(0,0,\ldots ,t_d,x)&+&H_d(0,0,\ldots  ,t_d, x, \frac{ \partial }{\partial x}u(0,0,
\ldots ,0,t_d,x))=0 \end{array}
\right.  \end{equation*}

The commutation of the $\widetilde\Phi_j^{t_j}$ implies that for any permutation $\sigma$ of
$\{1,....,n\}$, the Lagrangian manifold

$$L_{f;H_1,...,H_d}^{\sigma}=\bigcup_{(t_1, ... , t_q)\in{\mathbb  R}^q}\widetilde\Phi_{\sigma
(1)}^{t_{\sigma(1)}} ...
\widetilde\Phi_{\sigma (q)}^{t_{\sigma(q)}}(\Lambda_{ f;H_1,...,H_d})$$ coincides with $$L_{f;H_1,...,H_d}$$

Thus the associated function $u^\sigma_{f;H_1,...H_d}=u_{f;H_{\sigma(1)},...,H_{\sigma(d)}}$ coincides with
$u_{f,H_1,...,H_d}$. As a result, $u$ is a solution of any equation obtained by a replacing
$\frac{\partial}{\partial t_j}$ by
$\frac{\partial}{\partial t_{\sigma(j)}}$ and $H_j$ by $H_{\sigma(j)}$. Looking at the first equation, we see
that $$\frac{\partial}{\partial t_{\sigma(1)}}u(t_1,t_2,\ldots  ,t_d,x)+H_{\sigma(1)}(t_1,t_2,\ldots ,t_d, x,
\frac{\partial}{\partial  x}u(t_1,t_2, \ldots ,t_{d-1},t_d,x))=0$$

and since $\sigma(1)$ is arbitrary, we see that $u$ solves $(MHJ)$.

Our problem is thus solved for initial data and equation of class
$C^{1,1}$. 

However we are really interested in the $C^0$ case. This is the subject of the next section.

\subsection{The continuous case}

To deal with the $C^0$ case, we must first  extend the construction of  the Lagrangian submanifold
$$L_{f;H_1,...,H_d}$$
to  general  Hamiltonians and initial
conditions. Note that since the Hamiltonians do not commute anymore, order now matters.
 We shall consider for each permutation $\sigma$ of $\{1,...,d\}$,  the manifold
  $$L_{f;H_1,...,H_d}^{\sigma}$$
and as in Lemma 4.1, the manifold  $L_{f;H_1,...,H_d}^{\sigma}$ is Hamiltonianly isotopic to the zero
section, thus to  each such Lagrangian we may associate a unique  function $u^\sigma$.

If the $H_j$ are only $C^0$ and we consider a sequence $H_j^\nu$ of $C^2$ Hamiltonians converging to $H_j$. We shall then associate to these a sequence of Lagrangian submanifolds and a sequence of functions $u^{\sigma}_{\nu}$, which for $\sigma=Id$, is the
solution of the  set of equations obtained by replacing in \thetag{$\star$} the $H_j$ by $H_j^\nu$.

We thus want to prove that the conditions
$\lim_{\nu\to\infty}H_j^\nu = H_j$, and $\lim_{\nu\to\infty}\{H_j^\nu,H_k^\nu\}=0$,
 imply $$\lim_{\nu\to\infty} \vert u_{\sigma}^\nu-u_{\tau}^\nu
\vert =0$$ for any permutations $\sigma, \tau$. One should notice that this does not  say much about the
$L_{f;H_1^\nu,...,H_d^\nu}^{\sigma}$  and in particular we do not claim their convergence in the $C^0$ sense.

\begin{proposition} Assume $\lim_{\nu\to\infty}f^{\nu}=f$, $\lim_{\nu\to\infty}H_j^\nu = H_j$  and
$\lim_{\nu \to\infty}\{H_j^\nu,H_k^\nu\}=0$. Then  $u_{f;H_1,...,H_d}=\lim_{\nu\to\infty}  u_{\sigma}^\nu$ is
independent of $\sigma$ and is a solution of the multi-time Hamilton-Jacobi equation
\thetag{MHJ}. Moreover $u_{f;H_1,...,H_d}$ is a variational solution of each single equation.  If the $H_{j}$ are convex in $p$, then according to \cite{Zhukovskaia} variational and viscosity solutions do coincide  and thus $u_{f;H_1,...,H_d}$ is the viscosity solution of each single equation. \end{proposition}

The proof of this proposition will follow from a more general theorem, involving $c$-convergence. From now on, we shall only consider the case $d=2$, since the general case is in all respects similar.

\subsection{Capacity for manifolds with boundary}

In the sequel we shall need to define  $\gamma$ on $\mathscr L$ in the symplectic manifold $T^*([0,1]\times N)$ and  ${\mathscr D}_\omega(T^*([0,1]^2\times N)$ or more generally for the case of the cotangent bundle of a manifold $M$ with boundary.
Let $M$ be a manifold with boundary, and $L$ be a Lagrangian submanifold in $T^*N$ such that $\partial L =L \cap T_{\partial N}^*N$. Assuming moreover the intersection to be transverse,  we can consider the manifold $\tilde L \subset T^*(\tilde N)$ where $\tilde L, \tilde N$ are the doubles of $L, N$.
We shall write $\tilde N=N_-\cup N_+, \tilde L=L_-\cup L_+$ and we identify $L,N$ with $L_-,N_-$.
We also denote by $\tau$ be the involution of $\tilde N$ sending $N_-$ to $N_+$, and by $t$ the coordinate defining the collar of $\partial M$. By assumption, the reduction of $L$ by $\{t=t_0\}$, denoted by $L_{t_0}$ is a Lagrangian regular homotopy. We may modify this family near $\{t=0\}$, so that it is constant for $t$ close to $0$. Then near $\{t=0\}$, $L$ is uniquely defined by the real function $\tau(t,z)$ defined on $[0, \varepsilon]\times L_0$ such that $L=\{(t,\tau(t,z),z) \mid z \in L_t \}$.
Then $\tau(t,z)=\tau(z)$ near $t=0$, and we may thus extend $\tau$ and $L_t$ to $[- \varepsilon , \varepsilon ]\times L_0$, so that $\tau(-t,z)=\tau(t,z)$ and $L_{-t}=L_t$. This way we may extend
$L$ to a Lagrangian manifold $\tilde L \subset T^*\tilde M$, uniquely defined by $L$.

Notice that if $L$ is Hamiltonianly isotopic to the zero section,
by a Hamiltonian preserving $T_{\partial N}^*N$, then clearly the
same holds for $\tilde L$. Therefore $\tilde L$ has a GFQI, and we
can compute $\gamma (\tilde L)$. We define $\gamma
(L)=\gamma(\tilde L)$. Note also that if $S(x,\xi)$ is a GFQI for
$\tilde L$, then $ (\tau^*S)(x,\xi)=S(\tau(x),\xi)$ also generates
$\tilde L$.
In fact any $S$ generating $L$ over $T^*N$ has such an extension.

   \subsection{$S$-commuting Hamiltonians}
  \def \Ham {\mathscr {H}}
\def \DHam {\mathscr {HD}}

  Remember that we denoted by $\Ham (M)$ be the set of smooth time-dependent Hamiltonians on $ {\mathbb R} \times M$, and by $\DHam (M)$ the set of time one maps of such Hamiltonians, and defined ${\mathfrak H}_\gamma (M)$, the completion of $\Ham (M)$ for the metric $\gamma$, and $\mathfrak{HD}_\gamma (M)$ the completion of $\DHam (M)$ for the metric $\gamma$.

As we pointed out,  the ``time one flow'' map $\Ham (M)\longrightarrow \DHam (M)$
extends to a continuous map $\mathfrak {H}_\gamma (M) \longrightarrow
\mathfrak{HD}_\gamma (M)$, and in particular is well defined on the set of continuous maps $C^0( {\mathbb R} \times M) $.

  Let $\phi^s, \psi^t$ be the flows of the autonomous Hamiltonians  $H$ and $K$. Then we introduce the following commutation condition.
Consider $\rho^{(s,t)}=\phi^{s}\psi^{t}\phi^{-s}\psi^{-t}$ for $(s,t)\in [0,1]^2$, and the map

$$C(H,K): (s,t, \sigma,\tau, z) \to (s,t,\tilde\sigma (s,t,z), \tilde \tau(s,t,z), \rho^{(s,t)}(z))$$
from $T^*( [0,1] ^2\times N)$ into itself.
We may choose $\tilde\sigma, \tilde\tau$ such that $C(H,K)$ is a Hamiltonian diffeomorphism, using the fact that $\lambda \to \rho^{(s,\lambda t)}$ is an Hamiltonian isotopy starting from the identity. The  functions $\tilde\sigma, \tilde\tau$  are well defined up to the addition of $\frac{\partial c}{\partial s}(s,t)$ and $\frac{\partial c}{\partial t}(s,t)$  respectively, where $c$ is some $C^1$ function.

     If we are given  Hamiltonians $H(z), K(z)$ generating $\phi^s,\psi^t$, we can assume $s  \to \rho^{(s,t)}$ is generated by  $$ L_{t}(s,z)=H(z)-H(\psi^{-t} \phi^{-s}(z))$$

   We may also compute $$ \frac{ \partial }{ \partial t}L_{t}(s,z)= \{ H,K\}(\psi^{-t}\phi^{-s}(z))$$
and note that

 $$\tilde\tau_1(t_1,t_2,z)=\tau_1-L_{t_2}(t_1,\phi_1^{t_1}\phi_2^{t_2}\phi_1^{-t_1}\phi_2^{-t_2}(z))=H_1(\phi_1^{t_1}\phi_2^{t_2}\phi_1^{-t_1}\phi_2^{-t_2}(z))) - H_1(\phi_1^{-t_1}\phi_2^{-t_2}(z))$$

   $$ \tilde\tau_2(t_1,t_2,z)=\tau_2-t_1\int_0^1 \{H_1,H_2\}(\phi_1^{-\sigma t_1}\phi_2^{-t_2}(z))d\sigma
   $$

 Note that if either $M$ is compact and we impose the condition $\int_MH_j(z)\omega^n$
 or $H_j$ are compact supported, then $L_{t_2}(t_1,z)$ will satisfy the same condition, and this defines uniquely  $C(H,K)$ as the time one map of a Hamiltonian satisfying this normalization condition.

\begin{definition} \begin{enumerate}
\item
We shall say that $\phi^s, \psi^t$ $S$-commute up to $ \varepsilon $, if and only if $\gamma (C(H,K)) \leq \varepsilon $.
\item  Let $H,K$ be in $\mathfrak{HD}_\gamma (T^*N)$. We say that they $S$-commute, if and only if there are sequences $H_j, K_j$ such that $H_j,K_j$ converge to $H,K$ (for the metric $\gamma$) and  $\gamma (C(H_j,K_j))$ goes to zero.

\end{enumerate}
\end{definition}
We now have

 \begin{lemma} \begin{enumerate}
\item  If $\phi_1^t, \phi_2^t$ $S$-commute up to $ \varepsilon $ then $\gamma(\rho ^{(t_1,t_2)})\leq \varepsilon $
for all pairs $(t_1,t_2)\in [0,1]^2$. In particular $S$-commuting Hamiltonians $c$-commute.

\item If for all $(t_1,t_2)$ in $[0,1]^2$,  the Hofer norm of $\rho ^{(t_1,t_2)}$, is less than $ \varepsilon $, then $\phi_1^t, \phi_2^t$ $c$-commute up to $ \varepsilon $. In particular $C^0$-commuting Hamiltonians $S$-commute.
\end{enumerate}
 \end{lemma}
   \begin{proof} The first statement follows easily by symplectic 
    and the reduction inequality (see \cite{Viterbo-STAGGF}, page 705, prop. 5.1) that can be stated as follows: 
    
    Let $S(x,y,\xi)$ be a GFQI, and denote by $S_x$ the function $S_x(y,\xi)=S(x,y,\xi)$. 
    
    Then $$c(1,S)\leq \inf_x c(1_x,S_x) \leq sup_x c(\mu_x,S) \leq c(\mu,S)$$  
   
   The second statement follows from the proof of the inequality $$d(U) \leq 2\sup_x d(U_x)$$ in \cite{Viterbo-isop}, where $d$ is the displacement energy, i.e.
   
   $$d(U) = \inf \{ \vert H \vert _{C^0} \mid \phi^1_H(U)\cap  U=\emptyset \}$$  
  \end{proof}

\subsection{Different suspensions of Hamiltonian isotopies}

Let $H(t,z)$ be a time-dependent Hamiltonian on $T^*N$. We may associate to $H$ the Hamiltonians on $T^*( {\mathbb R} \times N)$ given by $\tilde H(t,\tau,z)=\tau+H(t,z)$ and $\widehat H(s;t,\tau,z)=sH(st,z)$, with respective flows

$$\Phi^s: (t,\tau,z) \to (t+s,\tau-H(t,\phi_t^{t+s}(z)), \phi_t^{t+s}(z))$$
$$\widehat \Phi^s: (t,\tau,z) \to (t,\tau-sH(st,\phi_0^{st}(z)), \phi_0^{st}(z))$$

Note that if $H$ is $C^0$ small, we have, denoting by $T_s$ the translation on the $t$ variable,   $\gamma(\Phi^sT_{-s})$ and $\gamma({\widehat \Phi}^s)$ are  small, in other words, $\Phi^s$ is $c$-close to $T_s$ and ${\widehat \Phi}^s$ is $c$-close to $Id$.

\section{Proof of theorem 1.3.}

The construction of variational solutions for a single equation can be rephrased  as follows.
Let $\Lambda_0$ be a Lagrangian submanifold of $T^*N$, and $H(t,z)$ be a time-dependent Hamiltonian. We consider $\tilde \Lambda_{0,H}=\{(0,-H(0,z),z) \mid z \in \Lambda_0\}$, and
$L= \bigcup_{t\in \R} \Phi^t(\Lambda_{0,H})$.

We can obtain $L$ as follows. Consider the map \begin{gather*} F: T^*(\R\times {\mathbb R} \times N) \to T^*(\R\times {\mathbb R} \times N)\\
(s,\sigma, \zeta)\mapsto (s,\tilde\sigma , \Phi^s(\zeta))
\end{gather*}
 and $L_{0,H}$ any Lagrangian submanifold in $T^*(\R \times N)$ such that $\tilde\Lambda_{0,H} \subset L_{0,H}$.
 Let  $t: T^*(\R \times N) \to \R$ be the projection on the first coordinate of $\R \times N$.

\begin{proposition} $L$ is the reduction of $ F(L_{0,H})$ by the coisotropic submanifold $\{s=t\}$.
\end{proposition}
 \begin{proof}
 Indeed, we look for the set of $(s,\tilde \sigma , \Phi^s(0,-H(0,z),z))$ such that $t=s$. But
 $\Phi^s (t,\tau,z)=(t+s,\tau +H(t,z)-H(t+s,\phi_t^{t+s}(z)), \phi_t^{t+s}(z))$. Intersecting with $t=s$ means that we consider the points   $\Phi^s(0,\tau,z)$ where $(0,\tau,z) \in L_{0,H}$ and by assumption this means that
 $(0,\tau, z) \in \tilde\Lambda_{0,H}$.
 \end{proof}

Now we will first  estimate $\gamma (L_{0,H})$ and then deduce the inequality $$\gamma(L)\leq \tilde\gamma (F)+ \gamma (L_{0,H})$$

\begin{lemma}\label{lemma-4.7}
Any isotopy $\Xi^s$ of $\{t=0\}$ in $T^*( {\mathbb R} \times N)$ can be generated by an arbitrarily small Hamiltonian preserving the hypersurfaces $\{t=t_0\}$.
\end{lemma}

\begin{proof}
   Indeed if we have  an isotopy of the type $(0,\tau, x, p) \mapsto  (0, \tau + \tau^s(x,p),x,p)$, it is induced by a Hamiltonian $H(t,x,p)$ and only depends on  $ \frac{ \partial H}{ \partial t}(0,x,p)= \frac{ \partial }{ \partial s}  \tau^s(x,p)$. The extension of $H$ to $T^*(\R \times N)$ can be arbitrarily small, for example we may take $H(s; t,\tau,x,p)=\chi(t)   \frac{ \partial }{ \partial s}  \tau^s(x,p)$, where $\dot\chi(0)=1$ and $\chi$ has arbitrarily small  $C^0$ norm.
  \end{proof}
\begin{proposition} For any positive $ \varepsilon $, we may choose $L_{0,H}$ so that we have
$$\gamma (L_{0,H})\leq \gamma (\Lambda_0) + \varepsilon $$
 \end{proposition}

 \begin{remark}
 Note that necessarily we have $\gamma (L_{0,H}) \geq \gamma ( \Lambda_0)$. This follows from the reduction inequality.
 \end{remark}
 \begin{proof}
 We notice that $L_{0,0}= \R \times \{0\}\times \Lambda_0$ satisfies obviously $\gamma (L_{0,0})=\gamma(\Lambda_0)$. Using the notation $L_{0,0}\cap \{t=0\}= \tilde \Lambda_{0,0}$,  the previous lemma allows us to find a Hamiltonian isotopy, $\Xi^s$ on $T^*( {\mathbb R} \times N)$ such that

 - $\Xi^s$ preserves $t=0$

 - on $t=0$, we have $\Xi^1(\tilde \Lambda_{0,0})=\tilde\Lambda_{0,H}$

 - $\Xi^s$ is generated by an arbitrarily $C^0$ small Hamiltonian.

Therefore $$\gamma(L_{0,H}) = \gamma(\Xi^1(\Lambda_{0,0})) \leq \tilde \gamma(\Xi ^1) + \gamma(\Lambda_{0,0}) \leq  \varepsilon + \gamma (\Lambda_{0,0})$$
 \end{proof}
\begin{Cor}
 We have $\gamma (L) \leq \tilde\gamma (F) + \gamma (\Lambda_0) $. Moreover if $L_1$ is the Lagrangian submanifold associated to $H_1, \Lambda_1$ and $L_2$ to $H_2,\Lambda_2$, we have

 $$\gamma (L_1-L_2)\leq \tilde\gamma (F_1\circ F_2^{-1})+ \gamma (\Lambda_1-\Lambda_2)$$

 \end{Cor}

  \begin{proof} The first inequality follows from the second by setting $$L_1=L, F_1=F, \Lambda_1=\Lambda, L_2=0_{\R\times N}, F_2=\Id, \Lambda_2=O_N$$
  Since  $L_1$ and $L_2$ are the reductions of $F_1(L_{1,H_1})$ and $F_2(L_{1,H_1})$, using the reduction inequality it will be enough to prove that
  $$\gamma (F_1(L_{1,H_1})-F_2(L_{2,H_2}))\leq \tilde\gamma (F_1\circ F_2^{-1})+ \gamma (\Lambda_1-\Lambda_2)$$
Since the map $F_2\circ F_1^{-1}$ sends $F_1(L_{1,H_1})$ to $F_2(L_{1,H_1})$, we have by the triangle inequality  and the symplectic invariance of $\gamma$,
 \begin{gather*}\gamma(F_1(L_{1,H_1})-F_2(L_{2,H_2}))\leq \gamma(F_1(L_{1,H_1})-F_2(L_{1,H_1}))+\gamma (F_2(L_{1,H_1})-F_2(L_{2,H_2})\leq  \\ \tilde\gamma(F_2\circ F_1^{-1})+ \gamma (L_{1,H_1}-L_{2,H_2})\end{gather*}
so that using invariance by reduction , we have
  $$\gamma (L_1-L_2)\leq \tilde\gamma(F_2\circ F_1^{-1}) + \gamma(L_{1,H_1}-L_{2,H_2})$$

  Finally we claim that for any positive $ \varepsilon$ we may find $L_{1,H_1}, L_{2,H_2}$ such that $\gamma (L_{1,H_1}-L_{2,H_2})\leq \gamma(\Lambda_1-\Lambda_2) + \varepsilon $. Indeed remember that $L_{1,H_1}, L_{2,H_2}$ are not uniquely defined, but just required to be Lagrangian submanifolds contained in $\tau+H_1(t,z)=0$ and $\tau+H_2(t,z)=0$ respectively, with reductions $L_1,L_2$. Since $\gamma (L_{1,0}-L_{2,0})=\gamma (\Lambda_1-\Lambda_2)$ and we saw in Lemma \ref{lemma-4.7} that
 we may go from $L_{1,0}$ to $L_{1,H_1}$ by applying a Hamiltonian with $C^0$ norm less than $ \varepsilon $, and similarly from $L_{2,0}$ to $L_{2,H_1}$
 we get that
 $ \gamma(L_{1,H_1}-L_{1,0})$ and $\gamma(L_{2,H_2}-L_{2,0})$ are  arbitrarily small.
 We may thus conclude
$$
  \gamma(L_{1,H_1}-L_{2,H_2})\leq \gamma(L_{1,H_1}-L_{1,0})+\gamma (L_{1,0}-L_{2,0})+\gamma (L_{2,0}-L_{2,H_2}) \leq 2 \varepsilon + \gamma (\Lambda_1-\Lambda_2)
  $$
Since $ \varepsilon$ is arbitrarily small, we get the announced inequality.

   \end{proof}
Our goal here is to prove the following theorem

 \begin{theorem} \label{strong-main-thm} Assume the  Hamiltonians $H_1(t_1,..,t_d,x,p),...,H_d(t_1,..,t_d,x,p)$ are Lipschitz, and {\bf
$S$-commute}. Then equation \thetag{MHJ} has a unique variational solution.  If the $H_j$ are convex in
$p$, then $u$ is a viscosity  solution of each individual equation.
 \end{theorem}

We limit ourselves to the case $d=2$.

We consider  the flows $\Phi_1^s, \Phi_2^s$ associated to  $\tau_1+H_1(t_1,t_2,z), \tau_2+H_2(t_1,t_2,z)$, and associated isotopies
 depending on $2$ parameters. More precisely  $\phi_1^{(t_1,t_1+s)}(t_2;z) $ is the flow of $(t_1,z)\to H(t_1,t_2,z)$ between times
  $t_1$ and $t_1+s$ ( $t_2$ being considered as fixed). Similarly
$\phi_2^{(t_2,t_2+s)}(t_1;z)$ is the flow associated to the Hamiltonian  $(t_2,z)\to H_2(t_1,t_2,z)$ between times $t_2$ and $t_2+s$.

Then $\Phi_1^s, \Phi_2^s$ are symplectomorphisms from $T^*({\mathbb R}^2 \times N)$ given by
$$\Phi_1^s (t_1,t_2,\tau_1,\tau_2,z)=(t_1+s,t_2,\tau_1+\hat\tau_1^s(t_1,t_2,z), \tau_2+\hat\tau_2^s(t_1,t_2,z),\phi_{1}^{(t_1,t_1+s)}(t_2;z))$$

$$\Phi_2^s (t_1,t_2,\tau_1,\tau_2,z)=(t_1+s,t_2,\tau_1+\check\tau_1^s(t_1,t_2,z), \tau_2+\check\tau_2^s(t_1,t_2,z),\phi_{2}^{(t_2,t_2+s)}(t_1;z))$$

Note that since $L_1=\bigcup_{t_1\in \R} \Phi_1^{t_1}(\tilde\Lambda_0)$ is not contained in the hypersurface $\tau_2+H_2(t_1,t_2,x,p)=0$, in order to get a Lagrangian submanifold after applying the flow $\Phi_2^s$ we must move $L_1$ in the $ \frac{\partial}{\partial \tau_2}$ direction, so that  its image is included in this hypersurface.

Setting $$\tau_2^{s}(z)=-\int_0^s\frac{ \partial H_1}{\partial t_2}(\sigma,0,\phi_1^{(0,\sigma)}(0;z))d\sigma$$ and
$$\widetilde L_1=\{(t_1,0,-H_1(t_1,0,\varphi^{(0,t_1)}(z)), -H_2(0,0,z)+\tau^{t_1}_2(z), \varphi_1^{(0,t_1)}(0;z))\mid z = (x, df(x))\in T^*N \}$$
 we may indeed check that on $\widetilde L_1$, $\tau_2+H_2(t_1,t_2,z)$ vanishes (note that $\widetilde \Lambda_1$ is contained in $\{t_2=0\}$).

We thus define a ``normalization map'', defined over $\{t_2=0\}$ :

$$V_{1,2}: (t_1,0,\tau_1, \tau_2,z) \mapsto (t_1,0,\tau_1,-H_2(t_1,0,z), z)$$

and similarly

$$V_{2,1}: (0,t_2,\tau_1, \tau_2,z) \mapsto (0,t_2,-H_1(t_1,0,z),\tau_2, z)$$

These maps are not diffeomorphisms, but their restriction to $\widetilde L_1$ (resp. $\widetilde L_2$) is one, and since it preserves the Liouville  form $(\tau_2dt_2+\tau_1dt_1+ \lambda)$, its restriction to $\widetilde L_1$ (resp. $\widetilde L_2$)  extends to a Hamiltonian isotopy $\tilde V_{1,2}$ (resp. $\tilde
V_{2,1}$) which according to lemma \ref{lemma-4.7} are generated by
arbitrarily $C^0$-small Hamiltonians \footnote{Note that $\widetilde V_{1,2}$ and $\widetilde V_{2,1}$ depend on $\widetilde L_1$ and $\widetilde L_2$ respectively.}.  We will use the maps

$$W_{1,2}: (s_1,s_2,\sigma_1,\sigma_2,\zeta) \mapsto (s_1,s_2,\hat\sigma_1, \hat\sigma_2, \Phi_2^{s_2}\tilde V_{1,2}\Phi^{-s_2})$$
and

 $$W_{2,1}: (s_1,s_2,\sigma_1,\sigma_2,\zeta) \mapsto (s_1,s_2,\check\sigma_1, \check\sigma_2, \Phi_1^{s_1}\tilde V_{2,1}\Phi_1^{-s_1})$$

 also generated by $C^0$-small Hamiltonians.

Now we claim that the following submanifolds are Lagrangian
 $$L_{1,2}=\bigcup_{(s_1,s_2)\in {\mathbb R} ^2} \Phi_1^{s_1}\tilde V_{2,1}\Phi_2^{s_2} ( \Lambda_0) $$
 $$L_{2,1}=\bigcup_{(s_1,s_2)\in {\mathbb R} ^2} \Phi_2^{s_2}\tilde V_{1,2}\Phi_1^{s_1} ( \Lambda_0) $$

It will later be useful to notice that  $$\Phi_2^{s_2}\tilde
V_{1,2}\Phi_1^{s_1}\Phi_2^{-s_2}\tilde V_{2,1}^{-1} \Phi_1^{-s_1}
=
  \left ( \Phi_2^{s_2}\tilde V_{1,2}\Phi_2^{-s_2}\right ) \circ  \left (\Phi_2^{s_2}\Phi_1^{s_1}\Phi_2^{-s_2} \Phi_1^{-s_1} \right ) \circ \left (\Phi_1^{s_1}\tilde V_{2,1}^{-1}\Phi_1^{-s_1}\right ) $$

 We then consider the map $C$ as above:

  \begin{gather*} T^*( {\mathbb R} ^2\times {\mathbb R} ^2\times N) \longrightarrow   T^*( {\mathbb R} ^2\times {\mathbb R} ^2\times N) \\ (s_1,s_2, \sigma_1, \sigma_2, \zeta) \longrightarrow  (s_1,s_2,\tilde\sigma_1, \tilde\sigma_2, \Phi_1^{s_1}\Phi_2^{s_2}\Phi_1^{-s_1}\Phi_2^{-s_2}(\zeta))
 \end{gather*}
 where $\zeta=(t_1,t_2,\tau_1,\tau_2,z)$,  and the map $\tilde C=W_{1,2}\circ C \circ W_{2,1}^{-1}$.

 We do not claim that $\tilde C$ sends $L_{1,2}$ to $L_{2,1}$, since these manifolds are in $T^*( {\mathbb R} ^2\times N)$ and not $T^*( {\mathbb R} ^2 \times {\mathbb R} ^2\times N)$.

However, let $\Delta$ be the coisotropic
 submanifold of $T^*( {\mathbb R} ^2\times {\mathbb R} ^2\times N)$ given by  $$\Delta=\{s_1=t_1,s_2=t_2\}$$

 \begin{proposition}\label{prop. 4.12}
 $$[{\tilde C}(0_{{\mathbb R}^2}\times L_{2,1})]_{\Delta} = L_{1,2}$$
  \end{proposition}
 \begin{proof}

 $$(0_{ {\mathbb R}^2}\times L_{1,2})= \{(s_1,s_2,0,0,\Phi_1^{t_1}\tilde V_{2,1}\Phi_2^{t_2}(\zeta))\mid \zeta \in \tilde \Lambda_{0}\}$$
 then
 \begin{gather*}\tilde C(0_{ {\mathbb R}^2}\times L_{2,1}) = \{(s_1,s_2,\tilde \sigma_1,\tilde \sigma_2,
 \Phi_2^{s_2}\tilde V_{1,2}\Phi_1^{s_1}\Phi_2^{-s_2}\tilde V_{2,1}^{-1} \Phi_1^{-s_1}(\zeta))\mid \zeta \in L_{2,1}\}=\\
  \{(s_1,s_2,\tilde \sigma_1,\tilde \sigma_2, (\Phi_2^{s_2}\tilde V_{1,2}\Phi_1^{s_1}\Phi_2^{-s_2}\tilde V_{2,1}^{-1}
   \Phi_1^{-s_1})\circ (\Phi_1^{t_1}\tilde V_{2,1}\Phi_2^{t_2})(\zeta))\mid \zeta \in \tilde \Lambda_0, (t_1,t_2) \in \R ^2 \}
  \end{gather*}

 Now if $t_1(\zeta)=t_2(\zeta)=0$ we have
 $$t_1((\Phi_2^{s_2}\tilde V_{1,2}\Phi_1^{s_1}\Phi_2^{-s_2}\tilde V_{2,1}^{-1} \Phi_1^{-s_1})\circ (\Phi_1^{t_1}
 \tilde V_{2,1}\Phi_2^{t_2})(\zeta))=t_1$$
 $$ t_2((\Phi_2^{s_2}\tilde V_{1,2}\Phi_1^{s_1}\Phi_2^{-s_2}\tilde V_{2,1}^{-1} \Phi_1^{-s_1})\circ
 (\Phi_1^{t_1}\tilde V_{2,1}\Phi_2^{t_2})(\zeta))=t_2$$

 Thus intersection of this with $\{s_1=t_1, s_2=t_2\}$ will be given by

 \begin{gather*}   \{(s_1,s_2,\tilde \sigma_1,\tilde \sigma_2,\Phi_2^{s_2}\tilde V_{1,2}\Phi_1^{s_1}\Phi_2^{-s_2}\tilde V_{2,1}^{-1} \Phi_1^{-s_1}\circ \Phi_1^{s_1}\tilde V_{2,1}\Phi_2^{s_2}(\zeta)\mid \zeta \in \tilde \Lambda_0\}=\\
  \{(s_1,s_2,\tilde \sigma_1,\tilde \sigma_2,\Phi_2^{s_2}\tilde V_{1,2}\Phi_1^{s_1}(\zeta)\mid \zeta \in \tilde \Lambda_0\}
 \end{gather*}

 and we must figure out its  projection according to the equivalence relation given by $$(s_1,s_2,\sigma_1,\sigma_2,t_1,t_2,\tau_1,\tau_2,z)\sim (s'_1,s'_2,\sigma'_1,\sigma'_2,t'_1,t'_2,\tau'_1,\tau'_2,z')$$
   if and only if
   $$s'_1=s_1=t_1=t'_1,s_2=s'_2=t_2=t'_2, \sigma_1+\tau_1=\sigma'_1+\tau'_1, \sigma_2+\tau_2=\sigma'_2+\tau'_2 $$

 Since the  map
    $$(s_1,s_2,\sigma_1,\sigma_2,t_1,t_2,\tau_1,\tau_2,z)\longrightarrow (s_1-t_1,s_2-t_2,\sigma_1,\sigma_2,t_1,t_2,\sigma_1+\tau_1,\sigma_2+\tau_2,z)$$
    is a symplectomorphism, the reduction of $X={\tilde C}(0_{{\mathbb R}^2}\times L_{2,1}) $ is given by

   $$X_\Delta=\{(s_1,s_2, \sigma_1+\tau_1,\sigma_2+\tau_2,z) \mid (s_1,s_2,\sigma_1,\sigma_2,s_1,s_2, \tau_1,\tau_2,z)\in X \}$$
   and this yields $L_{1,2}$ as announced.
 \end{proof}
  \begin{proposition}
  Let $\phi_1^t, \phi_2^t$ $S$-commute up to $ \varepsilon $, then we have $$\gamma (L_{1,2}-L_{2,1})\leq \varepsilon $$
  \end{proposition}

  \begin{proof}
  Indeed it is enough to prove that   $\tilde \gamma (\tilde C)$ is small by assumption ($\Phi^t, \Psi^s$ $S$-commute up to  $\varepsilon$). Since $\gamma( C)$ is small, and  $W_{1,2}, W_{2,1}^{-1}$  are generated by $C^0$ small Hamiltonians, $\gamma (W_{1,2}), \gamma (W_{2,1})$ are also small, and using the triangle inequality,
we get that $\tilde\gamma (\tilde C) $ is small.

Now if $\tilde\gamma(\tilde C) \leq \varepsilon $
  $$\gamma( {\tilde C}( 0_{\R ^2}\times L_{2,1})-  0_{\R ^2}\times L_{1,2}) \leq \varepsilon $$

Using the reduction inequality we obtain
  $$\gamma ({\tilde C}( 0_{\R ^2}\times L_{2,1})-  0_{\R ^2}\times L_{1,2}) _\Delta )\leq \varepsilon $$

and according to the above proposition \ref{prop. 4.12}, we get
$$\gamma( L_{2,1}-L_{1,2})\leq \varepsilon $$

  \end{proof}

\begin{proof}[Proof of  theorem \ref{main-thm}]
 Indeed we just proved
that if $\phi_1^t, \phi_2^t$ $S$-commute up to $\varepsilon$ on
$[0,T]$, we have $\gamma (L_{1,2}-L_{2,1}) \leq \varepsilon$ and
thus, by the reduction inequality we have

$$ \vert u_{1,2}(t,x)- u_{2,1}(t,x) \vert \leq \epsilon $$
on $[0,T]\times N$.

Now assume we have sequences $H_1^\nu,H_2^\nu$ such that $H_1^\nu$
$c$-converges to $H_1$, $H_2^\nu$ $c$-converges $H_2$, and
$H_1^\nu,H_2^\nu$ $S$-commute up to $\epsilon_\nu$ where
$\lim_{\nu \to 0}\epsilon_\nu=0$. Then we have

$$ \vert u_{1,2}^\nu(t,x)- u_{2,1}^\nu(t,x) \vert \leq \epsilon_\nu $$
on $[0,T]\times N$.

Now we know that $u_{1,2}^\nu$ converges to a variational solution
of
 \begin{gather*}\left \{
\begin{array}{ll}
\frac{\partial}{\partial t_1}\hat v
(t_1,t_2,x)+H_1(t_1,t_2,x,\frac{\partial}{\partial x}\hat
v(t_1,t_2,x)=0
\\
\hat v(0,t_2,x)=\hat w(t_2,x)
\end{array}\right .
 \end{gather*}

where $\hat w$ is the variational solution of the equation

\begin{gather*}\left \{
\begin{array}{ll}
\frac{\partial}{\partial t_2}\hat w
(t_2,x)+H_2(0,t_2,x,\frac{\partial}{\partial x}\hat w(t_2,x)=0
\\
\hat w(0,x)=f(x)
\end{array}\right .
 \end{gather*}

 on the other hand, $u_{2,1}^\nu$ converges to a variational solution
of

\begin{gather*}\left \{
\begin{array}{ll}
\frac{\partial}{\partial t_2}\check v
(t_1,t_2,x)+H_2(t_1,t_2,x,\frac{\partial}{\partial x}\check
v(t_1,t_2,x)=0
\\
\check v(t_1,0,x)=\check w(t_1,x)
\end{array}\right .
 \end{gather*}

where $\check w$ is the variational solution of the equation

\begin{gather*}\left \{
\begin{array}{ll}
\frac{\partial}{\partial t_1}\check w
(t_1,x)+H_1(t_1,0,x,\frac{\partial}{\partial x}\check w(t_1,x)=0
\\
\check w(0,x)=f(x)
\end{array}\right .
 \end{gather*}

Since we proved that $u_{1,2}^\nu(t,x)- u_{2,1}^\nu(t,x)$ goes to
zero, the sequences have a common limit,$v$, and this will be a
variational solution of both

$$\frac{\partial}{\partial t_2} v
(t_1,t_2,x)+H_2(t_1,t_2,x,\frac{\partial}{\partial x}
v(t_1,t_2,x)=0$$

and

$$\frac{\partial}{\partial t_1} v
(t_1,t_2,x)+H_1(t_1,t_2,x,\frac{\partial}{\partial x}
v(t_1,t_2,x)=0$$ and such that $v(0,0,x)=f(x)$.

In other words, $v$ is a variational solution of \thetag{MHJ}.
\end{proof}

 \begin{remark}
 Note that the initial condition on the function $u$ need not be given by ``functions''. It is enough to have a ``$C^0$- coisotropic submanifold'' $C \subset T^*( {\mathbb R} ^d \times N)$,  that  is the Hausdorff limit of  submanifolds $C_\nu$ such that  $\sigma_{TC_\nu}$ goes to zero for the $C^0$ topology, and
 $ T^*( {\mathbb R} ^d \times N)\to {\mathbb R} ^d$ given by $(t_1,...,t_d, \tau_1,...,\tau_d, x,p) \to (t_1,..., t_d)$ when restricted to the leafs of the foliation is uniformly proper. We let the reader adapt such a generalization into our framework. We also mention that the notion of $C^0$-Lagrangian submanifold has already been used by several authors, related to the fact that under various assumptions, we may guarantee that a smooth manifold that is  the  $C^0$-limit of Lagrangian submanifolds has to  be Lagrangian (see \cite{Laudenbach-Sikorav-IMRN},\cite{Sik3},\cite{Viterbo-IMRN}).
 \end{remark}

  \newcounter{appendix}  \setcounter{section}{\theappendix}

\def\thesection{\Alph{section}}

\section{Appendix: Proof of proposition 2.2}\label{Pf}
 \begin{lemma} \label{triangle-lag} Let $L_1,L_2,L_{3} \in {\mathscr L}$. Then $$\gamma(L_1-L_{3})\leq \gamma
(L_1-L_{2})+\gamma(L_2-L_{3})$$
\end{lemma}

 \begin{proof}
 Assume first  that $L_2=0_N$. Then we must prove
$$\gamma(L_1-L_3)\leq \gamma(L_1)+\gamma(L_3)
$$
  But this follows from corollary 2.8 and Proposition 3.3 on page 693 of \cite{Viterbo-STAGGF} according to which
 $$c(1, -S)=-c(\mu,S)$$
 and
 $$c(uv,S_1+S_2)\geq c(u,S_1)+c(v,S_2)$$
Indeed, we have $\gamma(L)=\gamma(-L)$ and

 $$c(1,S_1-S_3)\geq c(1,S_1)+c(1,-S_3)=c(1,S_1)-c(\mu,S_3)$$
 $$c(1,S_3-S_1)\geq  c(1,S_3)+c(1,-S_1)=c(1,S_3)-c(\mu,S_1)$$

using that $c(1,S_3-S_1)=-c(\mu,S_1-S_3)$ this can be rewritten as

 $$-c(1,S_1-S_3)\leq c(\mu,S_3)-c(1,S_1)$$
 $$c(\mu,S_1-S_3)\leq  c(\mu,S_1)-c(1,S_3)$$

and adding the inequalities we get

 $$\gamma(S_1-S_3)\leq \gamma(S_1)+\gamma (S_3).$$

  Then we claim that
  \begin{equation} \label{3.5++}
  \gamma(L-L')=\gamma(\phi(L)-\phi(L'))
  \end{equation}
  for any $\phi$ in $\DHam$. Indeed, according to proposition 3.5 page 695 of
\cite{Viterbo-STAGGF}, if
$L'=\rho^{-1}(0_N)$, we have

 \begin{gather*}\label{3.5+} \begin{array}{l}\gamma(\phi(L)-\phi (L'))=\gamma(\phi(L)-\phi
\rho^{-1}(0_N))= \\ \gamma((\phi \rho^{-1})^{-1}(\phi(L))=\gamma(\rho(L))\end{array}
  \end{gather*}
 Since the last expression does not depend on the choice of $\phi$, we get $$\gamma(\phi(L)-\phi
(L'))=\gamma(L-L')$$  Let now $\phi$ in $\DHam$ be such that
$\phi(L_2)=0_N$, and let us apply the above equality to the $L_j$:

 \begin{gather*} \gamma(L_1-L_3)=\gamma(\phi(L_1)-\phi(L_3))\leq \gamma(\phi(L_1))+\gamma(\phi(L_3)) = \\
\gamma (\phi(L_1)-\phi(L_2))+\gamma(\phi(L_2)-\phi(L_3))=\gamma(L_1-L_2)+\gamma(L_2-L_3)
\end{gather*}
 \end{proof}

We are now ready to prove proposition \ref{2.3} stating the main properties of  $\tilde\gamma$.
 
  \begin{proof}[Proof of proposition \ref{2.3}]
 
 \begin{enumerate}
\item Assume we have for all $L$ in $\mathscr L$, $\gamma (\psi(L)-L)=0$. According to \cite{Viterbo-STAGGF}, corollary 2.3, we have $\psi(L)=L$. Now since this holds for any $L$, we must have $\psi=\Id$.

\item Indeed the change of variable in ${\mathscr L}$ given by $L =\psi(L')$ shows that

$$\widetilde\gamma(\psi^{-1})=\sup\{\gamma(\psi^{-1}(L)-L) \mid L\in {{\mathscr L}} \}=
\sup\{\gamma(L'-\psi(L')) \mid L'\in {{\mathscr L}} \} $$

Since $\gamma(-S)=\gamma(S)$ we see that the last expression equals $\widetilde\gamma(\psi)$

\item  Again by change of variable, we have

 \begin{gather*} \sup\{\gamma(\psi\phi(L)-L) \mid L\in {{\mathscr L}} \} \leq  \\ \sup\{\gamma(\psi\phi(L)-\phi(L)) +
\gamma(\phi(L)-L) \mid L\in {{\mathscr L}} \} \leq \\ \sup\{\gamma(\psi(L')-L')\mid L'\in {{\mathscr L}} \} +
\sup\{\gamma(\phi(L)-L) \mid L\in {{\mathscr L}} \}
\end{gather*}  The second inequality follows from  lemma \ref{triangle-lag}

\item Indeed let us make the change of variable $L'=\phi^{-1}(L)$, we get
$$\widetilde\gamma(\phi\psi\phi^{-1})=\sup\{\gamma((\phi\psi\phi^{-1})(L)-L) \mid L\in {{\mathscr L}} \} =
\sup\{\gamma(\phi\psi(L')-\phi(L'))
\mid L'\in {{\mathscr L}} \}$$

Now using the equality (\ref{3.5++}), we get
$$\gamma(\phi\psi(L')-\phi(L'))=\gamma(\psi(L')-L')$$
 and thus $$\sup\{\gamma(\phi\psi(L')-\phi(L'))
\mid L'\in {{\mathscr L}} \}=\sup \{\gamma(\psi(L')-L')
\mid L'\in {{\mathscr L}} \}=\widetilde\gamma(\psi)$$
 \end{enumerate}
 \end{proof}

  \section{Appendix : Symplectic invariants and the main theorem in the non-compact case}\label{appendix-noncompact}

The aim of this appendix is to extend  our
construction and results to the case where $N$ is a non-compact manifold.
Of course, if all our Hamiltonians are compact supported, the
statements and proofs are straightforward adaptations of the compact
case.  We refer to \cite{Barles} and \cite{Fathi-Maderna} for the study, form a different point of view, of Hamilton-Jacobi equations in the non-compact setting.
\medskip

Let us first consider the problem of solving a single equation in our variational framework.
For this we start from the isotropic manifold $\Lambda_f$, and consider $L= \bigcup_{t\in \R} \Phi^t(\Lambda_{f,H})$.
For this to be defined, we need that the Hamiltonian flow is complete, or at least that its restriction to $\Lambda_{0,H}$ is complete \footnote{This means that $\Phi^s(z)$ is defined for all $z$ in $\Lambda_{f,H}$.}

\medskip

For a general Hamiltonian, we cannot hope to define variational
solutions unless the associated flow is complete, but this condition
will not suffice, as the following example shows.

\medskip

Consider the equation

$$
\left\lbrace
\begin{array}{ll}
& \frac{\partial u}{\partial t} - x^{2} - \frac{1}{4}   \vert \frac{\partial u}{\partial
  x} \vert ^{2} = 0\\
& u (0,x) = 0\ .
\end{array}
\right.
$$
The variational solution will be $u(t,x) = \tan(t) x^{2}$ which goes to
infinity as $t$ goes to $\frac{\pi}{2}$.

\medskip

This happens because the Lagrangian submanifold associated to the equation is given by
$$
L= \left\{ (t,\tau,x,p)| \quad \exists a \in \mathbb{R}, x = a \cos t, p
  =2 a \sin t , \tau = a^{2}\right\}
$$
and this being vertical for $t= \frac{\pi}{2}$, cannot contain the
graph of the differential of some function.

\medskip

We see that the trouble comes from the fact that the flow at time
$t=\frac{\pi}{2}$ brings over $\{ x = 0\}$ points which are arbitrarily
far. To avoid this phenomenon we introduce the notion of ``finite
propagation speed Hamiltonian''. We shall see that such a flow
can be suitably approximated by compact supported flows, and such
approximation can be used to define the variational solutions of the associated
Hamilton-Jacobi equation.

\medskip

We denote by $\pi : T^{*}N \to N$ the projection, by $\mathcal{W}$
some exhaustion $({W}_{\mu})_{\mu \in \mathbb{N}}$ of $T^{*}N$
by neighborhoods of the zero section, and by $(\Omega_{\mu})_{\mu\in
  \mathbb{N}}$ some exhaustion\footnote{By exhaustion of a set $X$, we mean a non-decreasing sequence $(Y_\mu)_{\mu \in \mathbb{N}}$ of subsets of $X$ such that $\bigcup_{\mu \in \mathbb{N}} Y_\mu=X$.} of $N$.

\medskip

We shall assume that for all $\mu$, and all compact subsets $K$ in
$N$, ${W}_{\mu} \cap \pi^{-1}(K)$, is compact.

\medskip

Our first task will be to determine when a sequence
$(\varphi_{\mu}^{t})$ of {\bf compact supported} Hamiltonian
isotopies {\bf approximates} a given Hamiltonian isotopy,
$\psi^{t}$, from the variational solution point of view. This will be
based on

\begin{proposition}\label{proposition-A1}
Let $\sigma^{t}$ be a compact supported Hamiltonian isotopy,
$L=\sigma^{1}(0_{N})$, $\rho^{t}$ another compact supported
Hamiltonian isotopy generated by $H(t,x,p)$, and set $\rho =\rho^{1}$.
Assume we have $|H(t,x,p)| \leq \varepsilon$ on
$$
\mathcal{U}((\rho^{t})_{t\in [0,1]}, L, \Omega) =
\left\{ (t, \rho^{t}(z))| t\in [0,1] \ , z\in L\ ,\ \rho(z)\in
  \pi^{-1}(\Omega) \right\}
$$
for some open set $\Omega$ in $N$.

Then for all $x$ in $\Omega$ we have
$$
| u_{\rho(L)}(x) - u_{L}(x)| \leq \varepsilon
$$
In particular for $\varepsilon = 0$, this means that if $H$ vanishes
on
$$
\bigcup_{t\in [0,1]} \{t\}\times \left\{ \rho^{t}(L\cap
  \rho^{-1}(\pi^{-1}(\Omega))\right\}
$$
we have
$$
u_{\rho(L)}(x) = u_{L}(x) \text{ for } x \text{ in }
\Omega .
$$
\end{proposition}

\begin{proof}
We first assume $|H| \leq \varepsilon$ everywhere. Then according to
Proposition 2.6 we have $\gamma (\rho(L)-L)\leq \varepsilon$ and the
inequality $|c(1_{x}, L') -c(1_{x},L)|\leq \gamma (L'-L)$ implies
$$
|u_{\rho(L)} (x) - u_{L} (x) |\leq \gamma (\rho (L) -L)|\leq
\varepsilon .
$$

Now to deal with the general case, assume $(H_{\tau})_{\tau \in
  [0,1]}$ is a smooth one-parameter family of compact supported Hamiltonians with $H_{0}= H$,
 and for all $\tau$ in $[0,1]$
$$
H_{\tau} = H_{0} \text{ on } \mathcal{U} ((\rho^{t}), L, \Omega).
$$

We claim that $c(1_x,\rho_\tau(L))=c(1_x,L)$.

Now notice that
$c(1_{x},\rho_{\tau}(L))$ belongs to the set $C_{\tau}(x,L)$ of
critical values of
$$(x,p,\xi) \longrightarrow \int^{1}_{0} \left[ p(s) \dot x(s) - H_{\tau} (s,
x (s) , p(s))\right] ds + S (x(0) ,\xi)$$
where $S$ is a generating function quadratic at infinity for $L$.

The set $C_\tau(x,L)$  is equal to the set of values of
$$
\int^{1}_{0} \left[ p_{\tau}(s) \dot x_{\tau} (s) - H_{\tau} (s,
x_{\tau} (s) , p_{\tau} (s))\right] ds + F (x_{\tau}(0) , p(0))
$$
where
$$
(x_{\tau}(s) , p_{\tau} (s)) = \rho^{s}_{\tau} (x_{\tau} (0) ,p_{\tau}
(0)) \leqno(a)
$$
$$
(x_{\tau} (0) , p_{\tau} (0))\in L \leqno(b)
$$
$$
x_{\tau} (1) = x \leqno(c)
$$
$$
dF = \lambda \text{ on } L\leqno(d)
$$
Since by assumption, $z\in L$ and $\rho(z)\in \pi^{-1}(\Omega)$ imply
$H_{\tau} = H_{0}$ on $\{ (s,\rho^{s}(z)), s \in [0,1]\}$, we conclude
$\rho^{s}_{\tau} (z) =\rho^{s}(z)$.

\medskip

Thus for $x$ in $\Omega$, the set $C_{\tau}(x,L)$ does not depend on
$\tau$. Since $C_{\tau}(x,L)$ is the set of critical values of the
action, it has measure zero, and using the fact that
$u_{\rho_{\tau}(L)} (x) = c(1_{x},\rho_{\tau}(L))$ depends
continuously on $\tau$, we conclude that it must be independent of
$\tau$ for $x$ in $\Omega$.

>From the above argument follows  that
$$
u_{\rho_{\tau}(L)} (x) = u_{\rho(L)} (x) \text{ for all } x \text{ in }
\Omega\ .
$$
Now if $\vert H_{0}\vert \leq \varepsilon$ on $\mathcal{U}((\rho^{t})_{t\in[0,1]},L,\Omega)$
we may find a linear homotopy
$$
H_{\tau} (s,z) = (1-\tau \chi (s,z)) H_{0}(s,z)
$$
where $\chi(s,z)=0$ for $(s,z)\in \mathcal{U} ((\rho^{t}), L,\Omega)$
and $\chi =1$ outside a neighborhood of $\mathcal{U}((\rho^{t}),L,\Omega)$.

\medskip

By suitably choosing $\chi$,we may assume the following inequality holds $$|H_{1}| \leq
\varepsilon +\delta$$ with $\delta > 0$ arbitrarily small. Since
$H_{\tau}=H_{0}$ on $\mathcal{U}((\rho^{t}),L,\Omega)$, we have
$$u_{\rho_{1}(L)}(x) = u_{\rho(L)}(x)$$ for $x$ in $\Omega$. Since
$|H_{1}|\leq \varepsilon + \delta$ everywhere, we may apply the first part of the proof and obtain
$$
|u_{\rho^{1}(L)}(x) -u_{L}(x) |\leq \varepsilon +\delta
$$
and thus
$$
|u_{\rho(L)}(x) - u_{L}(x)|\leq \varepsilon + \delta\ .
$$
Finally, $\delta$ being arbitrarily small, we conclude that
$$
|u_{\rho(L)} (x) - u_{L}(x) | \leq \varepsilon
$$
\end{proof}
Our next task is to associate to the non compact supported Hamiltonian isotopy
$\psi^{t}$ a sequence $(\varphi^{t}_{\mu})_{\mu\in\mathbb{N}}$ of compact supported ones,
such that $c(1_x,\phi_\mu^t(L))$ converges -or rather its restriction to any compact set stabilizes- for fixed $x,t$, as
$\mu$ goes to infinity.

\medskip

We denote by $H_{\mu}(t,z)$ and $H(t,z)$ the Hamiltonians
generating $\varphi_{\mu}^{t}$ and $\psi^{t}$, and  by  $K_{\mu,\nu} (t,z)$ the Hamiltonian generating $\rho_{\mu
,\nu}^{t} = \varphi^{t}_{\nu}\varphi^{-t}_{\mu}$.

\medskip

The equality $\varphi^{t}_{\mu}= \varphi^{t}_{\nu} =\psi^{t}$ on
$\mathcal{U}$ is equivalent to $H_{\mu}(t,z) = H_{\nu}(t,z) =
H(t,z)$ on $\psi^{s}(\mathcal{U})$ for all $s\in [0,T]$.

Thus if $\varphi^{s}_{\mu} = \varphi^{s}_{\nu} = \psi^{s}$ on
$$
{W}_{\mu}\cap (\psi^{t})^{-1} (\pi^{-1} (\Omega_{\mu}))
$$
we have.
$$
\varphi_{\nu}^{s} \varphi_{\mu}^{-s} = id \text{ on }
\psi^{s}({W}_{\mu}\cap (\psi^{t})^{-1} \pi^{-1}(\Omega_{\mu}))
$$
hence
$$
K_{\mu,\nu}(s,z) = 0 \text{ on } \psi^{s}({W}_{\mu}\cap (\psi^{t})^{-1}
\pi^{-1} (\Omega_{\mu}))
$$
so that, according to proposition \ref{proposition-A1},
$$
u_{\varphi^{t}_{\mu}(L)} (x) = u_{\varphi_{\nu}^{t}(L)}(x) \text{
for } x \text{ in } \pi^{-1}(\Omega_{\mu})
\; t\in [0,T] $$

This motivates the
\begin{definition}
Let $(\psi^{t})_{t\in[0,T]}$ be a Hamiltonian isotopy. A sequence
of {\bf compact supported} Hamiltonian isotopies,
$(\varphi_{\mu}^{t})_{\mu\in\mathbb{N}}$ is an {\bf exhausting
sequence} associated to $\mathcal{W}$ and
$(\Omega_{\mu})_{\mu\in\mathbb{N}}$ if for all $\mu$, we have
$$
\varphi^{t}_{\mu} = \psi^{t} \text{ on } \bigcup_{s\in [0,T]}
{W}_{\mu} \cap (\psi^{s})^{-1}(\pi^{-1}(\Omega_{\mu}))
$$
\end{definition}

The above argument implies

\begin{proposition}\label{prop-A3}
Let $(\varphi_{\mu}^{t})$ be an exhausting sequence for
$(\psi^{t})_{t\in [0,T]}$ associated to $\mathcal{W}$ and
$(\Omega_{\mu})_{\mu\in\mathbb{N}}$. Then for any $L$,  image of the zero section by a compact supported Hamiltonian isotopy and contained in
some ${W}_{\mu}$, and, for all $\alpha ,\beta \geq \mu$
and $t\in [0,T]$
$$
u_{\varphi^{t}_{\alpha}(L)}(x) = u_{\varphi^{t}_{\beta}(L)}(x)
\quad \forall x\in \Omega_\mu
$$
\end{proposition}
We may therefore set

\begin{definition}
Given $(\psi^{t})_{t\in [0,T]}$ and an exhausting sequence $(\phi^t_\mu)$
associated to $\mathcal{W} , (\Omega_{\mu})_{\mu\in \mathbb{N}}$
we define $u_{\psi^{t}(L)}(x)$ to be the common value of the
$u_{\varphi^{t}_{\mu}(L)}(x)$ for $\mu$ large enough.
Since $u_{\psi^t(L)}$ is a solution of the Hamilton-Jacobi equation, we call it the {\bf variational solution associated to the exhausting sequence.} \end{definition}

We now need to solve two difficulties.

\medskip

First, it is not clear, given $(\psi^{t})_{t\in [0,T]}$,
$\mathcal{W}\ ,\ (\Omega_{\mu})_{\mu\in\mathbb{N}}$ whether we may
find an associated exhausting sequence. Second, one could wonder
whether $u_{\psi^{t}(L)}(x)$ will depend on the choice of the
exhausting sequence. But this will not be  the case, since if $\check\phi^t_\mu, \hat\phi_\mu^t$ are exhausting sequences, we get  a new exhausting sequence by intertwining them. Remember that we assume $W_\mu \cap
\Omega_\mu$ is compact for all $\mu$.

\begin{definition} \begin{enumerate} \item
 We shall say that $(\psi^{t})_{t\in[0,T]}$ has {\bf finite propagation
speed} with respect to $\mathcal{W}$ if and only if
 for each $\mu$, there exists $\nu$ such that
$$
\bigcup_{t\in[0,T]} \psi^{t}({W}_{\mu})\subset W_\nu
$$

\item  $\mathcal W$ is a {\bf convex exhaustion} if the fibers of $W_\mu$ are  convex, that is for all $x\in N$  $W_\mu\cap \pi^{-1}(x)$ is convex in $T^*_xN$.
\end{enumerate}
\end{definition}

Note that the composition of finite propagation speed Hamiltonians has
finite propagation speed (with respect to the same family
$\mathcal{W}$). However this is not necessarily true for its inverse
$((\psi^{t})^{-1})_{t\in [0,T]}$.

\medskip

The set of Hamiltonians isotopies having, together with their inverse, finite
propagation speed, is a group for the composition law $$(\phi^t) \circ (\psi^t) = (\phi^t\psi^t)$$

Our last result for this section is

\begin{proposition}
Let $\mathcal W$ be an exhaustion.
If $(\psi^{t})_{t\in [0,T]}$ has finite propagation speed with
respect to $\mathcal{W}$, then, for any exhaustion $(\Omega_\mu)_{\mu \in {\mathbb N}}$ of $N$, there exists an exhausting sequence
$(\varphi_{\mu}^{t})$ associated to  $\mathcal{W}$.
\end{proposition}

The proof makes use of the  obvious lemma:

\begin{lemma}

Let $K\subset U$ be a compact subset and $(\psi^{t})_{t\in[0,1]}$ a
Hamiltonian isotopy such that
$$
V = \bigcup_{t\in [0,1]} \psi^{t}(K) \subset U
$$
Then there is a Hamiltonian isotopy, $\varphi^{t}$, supported near $V$
such that
$$
\varphi^{t} = \psi^{t} \text{ on } K \text{ for all } t \text{ in
} [0,1]
$$
\end{lemma}

\begin{proof}
Indeed, if $\chi$ is a function equal to 1 on $V$, vanishing outside a
neighborhood of $U$, and $H(t,z)$ generates $\psi^{t}$, then $\chi(z)
H(t,z)$ generates $\varphi^{t}$.
\end{proof}

\begin{proof}
We may now prove the proposition.

Since for all $\mu$ there exists $\nu$ such that
$$\psi^s(W_\mu) \subset W_\nu$$

for all $s$ in $[0,T]$, we have

$$\psi^s(W_\mu) \cap \pi^{-1}(\Omega_\mu) \subset W_\nu \cap \pi^{-1}(\Omega_\mu)\subset
W_\nu \cap \pi^{-1}(\Omega_\nu) $$

Thus

\begin{gather*}
\bigcup_{s\in [0,T]}W_\mu \cap (\psi^s)^{-1}(\pi^{-1}(\Omega_\mu))
=\bigcup_{s\in [0,T]}(\psi^s)^{-1}(\psi^s(W_\mu) \cap \pi^{-1}(\Omega_\mu)) \subset
\bigcup_{s\in [0,T]}(\psi^s)^{-1}(K_\nu)\end{gather*}

We denote by $\tilde K_{\nu,T}$ this last set, and note that it is
compact. Therefore $$\bigcup_{s\in [0,T]}\psi^s(\tilde K_{\nu,T})$$ and we may find an isotopy $\phi^t_\mu$, compact supported, such that
$\phi_\mu^t=\psi^t$ on $\tilde K_{\nu,T}$.
 \end{proof}

Now this yields solutions of the evolution Hamilton-Jacobi equations  in case the initial condition $f$ has compact support. Of course formally the general case can be reduced to this one, since solving
 \begin{gather*}\left\{ \begin{array}{ll} \frac{\partial u}{\partial t}(t,x)+H(t,x,\frac{
\partial u}{\partial x}(t,x))=0 \\ u(0,x)=f(x) \end{array}\right .\end{gather*}

is equivalent to solving

\begin{gather*} \left\{\begin{array}{ll}
\frac{\partial v}{\partial t}(t,x)+H(t,x,\frac{
\partial v}{\partial x}(t,x)+df(x))=0 \\ v(0,x)=0
\end{array}\right. \end{gather*}

 by the change of function $v(t,x)=u(t,x)-f(x)$ (and this change of function preserves variational solutions).

 However we would like to have a criterion for $H$ and an independent one for $f$  in order to isolate  the difficulties.

  \begin{proposition}
  Let $\mathcal W$ be a convex exhaustion.  Let $ C_{\mathcal{W}}(N)$ be the union over all sets $W_\mu$ of the exhaustion, of the $C^0$ closure
   of the set of $C^1$ functions with graph in $W_\mu$.
   Let $ \mathcal { H}^1_{ \mathcal W}( {\mathbb R} \times T^*N)$ be the set of  $C^{1,1}$ Hamiltonians having
   finite propagation speed with respect to $\mathcal{W}$.
   Let  $f$ such that the graph of $df$ is in $W_\mu$,  the function $u_{\psi^{t}(L)}$ does not depend on the
choice of the exhausting sequence. Moreover it is a solution of the Hamilton-Jacobi equation
\begin{gather*}\tag{HJ} \left\{\begin{array}{ll}
\frac{\partial u}{\partial t}(t,x)+H(t,x,\frac{
\partial u}{\partial x}(t,x))=0 \\ u(0,x)=f(x)
\end{array}\right . \end{gather*}
 \end{proposition}

\begin{remark}
Let $\psi^t$ be the flow of $H$ and $\Lambda_0\subset T^*N$ be a Lagrangian submanifold.

We first notice that the manifold we  constructed, $L$ in $T^*(\R\times N)$ giving the solution of \thetag{HJ} and the family $L_t=\psi^t(L)$ give essentially the same information. Indeed, if $u(t,x)=c(1_{(t,x)},L)$ and $v(t,x)=c(1_x,L_t)$ we have $u(t,x)-v(t,x)=d(t)$ for some continuous function $d$.
Given $v$, we can easily determine $d$ or $u$, since $ \frac{\partial}{\partial x} v(t,x)= \frac{\partial}{\partial x} u(t,x)$, hence
$$ \frac{\partial}{\partial x} u(t,x) = -H(t,x, \frac{\partial}{\partial x}v(t,x))$$
hence
$$ u(t,x)=f(x)+\int_0^tH(s,x, \frac{\partial}{\partial x} v(s,x))ds$$
\end{remark}

\begin{proof}  Clearly if $\Lambda_0$ is the graph of $df$, we wish to replace it by a Lagrangian manifold $\Lambda_\mu$, such that there exists an isotopy $\Lambda_\mu^s$ from $\Lambda_0$ to $\Lambda_\mu$, such that $$ \phi_\mu^t(\Lambda_\mu^s)=\psi^t(\Lambda_\mu^s)=\psi^t(\Lambda_0)$$ over $\pi^{-1}(\Omega_\mu)$ for all $t \in [0,T]$.

 But $\phi_\mu^t=\psi^t$ on $W_\mu\cap (\psi^s)^{-1}\pi^{-1}(\Omega_\mu)$
  so if $\Lambda_\mu^s= \Lambda_0$ over $\pi^{-1}(A_\mu)$, where $A_\mu$ will be determined later, we have

 $$\phi_\mu^t(\Lambda_\mu^s)=\psi^t(\Lambda_0)$$ over $$\psi^t(W_\mu\cap (\psi^s)^{-1}\pi^{-1}(\Omega_\mu)\cap \pi^{-1}(A_\mu))$$

By compactness of $$W_\mu\cap (\psi^s)^{-1}\pi^{-1}(\Omega_\mu)$$ we may find $A_\mu$ such that $W_\mu\cap (\psi^s)^{-1}\pi^{-1}(\Omega_\mu)\subset \pi^{-1}(A_\mu)$,
and therefore according to proposition \ref{prop-A3}
$$u_{\phi^t_\mu(\Lambda_\mu^1)}=u_{\psi^t(\Lambda_\mu^1)}$$ on $\Omega_\mu$.

Now if $\Lambda_0$ is the graph of $df$, we may take for $\Lambda_\mu^s$ the graph of $df_\mu^s$ where $$f_\mu^s(x)=(1-s\chi_\mu(x))f(x)$$

and $$\left\{ \begin{array}{ll} \chi_\mu(x)=0 \;\text{on} \; A_\mu \\ \chi_\mu(x)=1 \;\text{near infinity}\; \end{array}\right .$$

Then $\Lambda^s_\mu$ will coincide with $\Lambda_0$ over $A_\mu$. We still need to prove that $\Lambda _\mu^s$ is in $W_\mu$ for all $s$ in $[0,1]$. Now this coincides with the graph of $(1-s\chi_\mu(x))df(x)$ outside a compact set. Since the graph of $df$ is in $W_\mu$, this will be also in some  $W_\nu$ for some $\nu$, provided $W_\nu$ is fiber convex.  

\end{proof}
 \begin{remark} In practice most $W_\mu$ are fiber convex, and if this is not the case, we can always convexify them. However, it is not clear that a $\psi^t$  having finite propagation speed with respect to $\mathcal {W}$ will have also finite propagation speed with respect to the convexified exhaustion $\mathcal {W} '$, unless  the convexification of $W_\mu$, is contained in some $W_\nu$.
 \end{remark}

\begin{example}
  \begin{enumerate}
  \item The Barles-Tourin condition defined  for an autonomous Hamiltonian $H(x,p)$ on $T^*( {\mathbb R} ^n)$ can
be rewritten as

 \begin{itemize} \item $$ \vert D_pH(x,p) \vert \leq f( \vert p \vert )(1+ \vert x \vert )$$
\item  $$ g_-( \vert p \vert )\leq \vert H(x,p) \vert \leq g_+(
\vert p \vert )$$
 where $\lim_{\xi\to\infty} g_{\pm}(\xi)=+\infty$
 \end{itemize}

 This implies finite propagation speed for the exhaustion by the  $U_r=\{(x,p) \mid \vert p \vert \leq r \}$,
since the assumptions and conservation of energy imply that the
image of $U_r=\{(x,p) \mid \vert p \vert \leq r \}$
 by the flow remains in  some $U_{r'}$. In such a region, the first inequality implies, by Gronwall's lemma,
finite propagation speed.
 \item Consider the $U_{(a,b)}=\{(x,p) \mid \vert p \vert \leq ad(x,x_0)+b $, and assume
  \begin{itemize}
\item  $$ \vert D_xH(t,x,p) \vert \leq C( \vert p \vert )$$ \item
$$ \vert D_pH(t,x,p) \vert \leq C(1+ d( x,x_0) )$$
 \end{itemize}
  Then the first assumption, together with Gronwall's lemma implies that  the flow sends $U_{a,b}$ to
  some $U_{a',b'}$.  \end{enumerate}
 \end{example}

The proposition implies that if $\mathcal W$ is an exhaustion, and $f$ a function on $N$ such that the graph of $df$ is contained in $ W_\mu$ for some $\mu$, and if $H(t,x,p)$ generates a flow having finite propagation speed with respect to the $\mathcal W$, then we may find a variational solution of
 \begin{gather*}
 \frac{\partial u}{\partial t}(t,x)+H(t,x,\frac{
\partial u}{\partial x}(t,x))=0 \\ u(0,x)=f(x)  \end{gather*}

by considering solutions of

 \begin{gather*}
\left \{ \begin{array}{ll} \frac{\partial u_\mu}{\partial t}(t,x)+H_\mu(t,x,\frac{
\partial }{\partial x}u_\mu(t,x))=0 \\ u_\mu(0,x)=f(x) \end{array}\right . \end{gather*}

where $H_\mu$ is a compact supported Hamiltonian generating an exhausting sequence $(\phi_\mu)_{\mu\in {\mathbb N}}$ for the flow $\psi^t$ of $H$.

  It follows from the above that $u(t,x)=\lim_{\mu\to\infty}u_\mu(t,x)$ defines a function independent of the choice of the exhausting sequence. In the regular case (i.e. for smooth data), this will of course be a solution of the equation almost everywhere.

 We may of course combine theorem \ref{strong-main-thm} with  the above proposition to get the following:

 \begin{theorem} \label{super-thm}
  Let $\mathcal W$ denote an exhausting sequence of neighborhoods of $0_N$.
   Let $ \mathfrak{C}_{\mathcal{W}}(N)$ be the union over all sets $W_\mu$ of the exhaustion, of the $\gamma$ closures
   of the set of $C^1$ functions with graph in $W_\mu$.
   Let $ \mathcal { H}^1_{ \mathcal W}( {\mathbb R} \times T^*N)$ be the set of  $C^{1,1}$ Hamiltonians having
   finite propagation speed with respect to $\mathcal{W}$.
Consider the closure $ \mathfrak{ H}_{\gamma, \mathcal W}( {\mathbb R}\times  T^*N)$ of $ \mathcal { H}^1_{ \mathcal
W}( {\mathbb R} \times T^*N)$ in $ \mathfrak{H}_\gamma( {\mathbb R} \times T^*N)$.
 Finally we assume that $$(f, H_1, ...., H_d) \in \mathfrak{C}(N)\times \mathfrak{H}_{\gamma, \mathcal W}( {\mathbb R} \times T^*N)^d$$
 then there is a solution of \thetag{MHJ} with data $(f,H_1,..., H_d)$.

 \end{theorem}

 \section{Appendix : equations involving the unknown function}\label{Appendix-B}

 For equations of the type $$\frac{\partial u}{\partial t_j}+H_j(t_1, ... , t_d,x,\frac{\partial}{\partial x}u(x),u(x))=0$$

 the same proofs hold, but we have to replace the cotangent space $T^*M$ which is symplectic by  the jet space $J^1(M)$, a contact manifold, on which we use the contact form $\alpha =dz-pdx$. We also replace  Lagrangian submanifolds, by Legendrian ones, which are  manifolds of dimension $\dim (M)$ on which $\alpha$ vanishes. It was noticed by Chekanov (around 1986, but published in \cite{Chekanov} much later), that the theory of generating functions holds also in the Legendrian case, and that existence is invariant by Legendrian isotopy. 

 Since the notion of generating function applies {\bf  ``even better''} to a Legendrian submanifold of $J^1(M)$, and according to
 \cite{Theret-these}, we also have uniqueness of generating functions, all the constructions of \cite{Viterbo-STAGGF} apply. Similarly, since the Hamiltonian  flow of $\tau + H(t,x,u,p)$ is also well-defined, the constructions from \cite{Viterbo-Ottolenghi} also carry through.

 The details of proofs are left to the reader. However, we will make explicit the commutativity condition for $H_1(t_1,t_2,x,p,u)$ and  $H_2(t_1,t_2,x,p,u)$.

 This is easy even in a general contact manifold. Indeed, define the Reeb vector field $R$, by the equalities
  \begin{gather*} i_{R}\alpha =1 \\ i_{ R} d\alpha = 0 \end{gather*}
 the  associated Hamiltonian vector field must satisfy
 \begin{gather*} i_{\tilde X_H}\alpha =-H \\ i_{\tilde X_H} d\alpha = dH+i_R(d H) \alpha  \end{gather*}
 (note that for $H=-1$ we recover $R$)
 which in our case, since  $ \alpha=du-pdx, R=\frac{\partial }{\partial u}$ translates to

 $${\tilde X}_H(x,p,u)=(pH_p-H)\frac{\partial }{\partial u} +H_p\frac{\partial }{\partial x} -(H_x+pH_u)\frac{\partial }{\partial p}$$

 We then have  $L_{{\tilde X}_H}\alpha =H_u \alpha$ and this shows that $\tilde X_H$ preserves the contact distribution $\alpha=0$. Since the commutator of such vector fields must also preserve the contact distribution, we may define

 \begin{definition}
 Let  $H_1(x,p,u)$ and  $H_2(x,p,u)$ be two functions on $J^1(M)$.  We set
  \begin{gather*} [H_1,H_2](x,p,u)=\\ \{H_1,H_2\}_s(x,p,u)+\\ \frac{\partial H_2}{\partial u}(x,p,u)\left (p\frac{\partial H_1}{\partial p}(x,p,u)-H_1(x,p,u)\right )- \frac{\partial H_1}{\partial u}(x,p,u) \left (p\frac{\partial H_2}{\partial p}(x,p,u)-H_2(x,p,u)\right ) \end{gather*}
 where $$ \{H_1,H_2\}_s = \frac{\partial H_1}{\partial p}\frac{\partial H_2}{\partial x} -\frac{\partial H_1}{\partial x}\frac{\partial H_2}{\partial p}$$
 \end{definition}

 \begin{proposition}
$${ \widetilde X}_{[H_1,H_2]}=[\widetilde X_{ H_1}, \widetilde X_ { H_2}]$$

 \end{proposition}

Now in the multi-time dependent case, for $H_1(t_1,t_2,x,p,u), H_2(t_1,t_2,x,p,u)$, we can consider $$K_1(t_1,t_2,\tau_1,\tau_2,x,p,u)=\tau_1+
 H_1(t_1,t_2,x,p,u), K_2(t_1,t_2,\tau_1,\tau_2,x,p,u)=\tau_2+ H_2(t_1,t_2,x,p,u)$$
as functions on $J^1({\mathbb R}^2\times M)$  and set

   \begin{definition}
   Let $H_1,H_2$ as above. We define
   $$[[H_1,H_2]]= [K_1,K_2]$$
   \end{definition}

The following is the natural extension of Theorem 1.2 :
  \begin{proposition}
Assume we have sequences $H_\nu, K_\nu$ such that $H_\nu {\overset {C^0} \longrightarrow } H, K_\nu {\overset {C^0} \longrightarrow } K, [H_\nu, K_\nu] {\overset {C^0} \longrightarrow } 0$. Then
$[H,K]=0$.
  \end{proposition}

\section { Appendix : Coisotropic submanifolds of symplectic manifolds}
\medskip
\par\noindent
Let $(P,\sigma )$ be a symplectic manifold, and $C$ a coisotropic submanifold, that is
$$
(T_z C)^{\perp}\subset T_z C\qquad \forall z\in C
$$
where
$$
(T_z C)^{\perp}:= \{ v\in T_z P:\ \ \sigma (v,u)=0\ \  \forall u\in T_z C   ,\}
$$
here $\sigma$ is a closed and non
degenerated 2-form, e.g. if
$P=T^*M$ , then $\sigma =d\theta$,  $\theta=\sum_{i}p_idx^i$ is the  Liouville 1-form.
Let $H:T^*M\to {\mathbb R}$, then the Hamiltonian vector field  $X_H$ related to $H$
is so defined:
$$
\sigma (X_H, \cdot )=-dH,
$$
and the Poisson brackets of two Hamiltonian functions $H$ e $K$ are
$$
\{K,H\}(z)=dK(z)\cdot X_H(z)\qquad \quad (z=(x,p)).
$$
Let suppose that
$C$ is defined by zero's of functions\footnote{it is
equivalent to ask that the normal bundle of $C$, with respect to some
Riemannian metric, is trivial (\cite{Guillemin-Pollack}  p. 77.) }, $H^a :P\to
{\mathbb R}$,
 $a=1,...,{\rm codim}\,C$:
$$
C:\ \ H^a =0\qquad {\rm rk}\,dH^a\big|_{H^a =0}={\rm max}=
{\rm codim}\,C
$$
Thus: $u\in T_z C \iff dH^a\, u =0$. Standard homomorphism theorem
between linear spaces gives existence (and uniqueness, by the max
rank) of the Lagrange multipliers $\lambda_a=\lambda_a (z,v)$:
$$
\forall v\in (T_zC)^{\perp}:\quad
\big[dH^au=0\Rightarrow \sigma (v,u)=0\big]\Rightarrow \omega (v,\cdot
)=\sum_a\lambda_a dH^a
$$
The linear independence of the forms $dH^a,\  a=1,..., {\rm codim}C$, implies that
the following Hamiltonian  vector fields $\Big( X_a\Big)_{a=1,..., {\rm codim}C}$,
$$
\sigma (X_a,\cdot
)=- dH^a,
$$
span $(T_z C)^{\perp}$. Since $(T_z C)^{\perp}\subset T_z C$, the vectors
$X_a (z)$ are in
$T_zC$, and
their integral paths, the {\em characteristics} of $X_a$, are in $C$,
hence
$$
\sigma (X_a,X_b)=0.
$$
>From the definition of Poisson brackets $\{\cdot ,\cdot \}$,
$$
\sigma (X_a,X_b)=\{H_a, H_b\},
$$
and from the well known Lie algebra morphism ($[\cdot , \cdot ]$: Lie brackets)
$$
[X_a,X_b]=X_{ \{H_a,H_b\} },
$$
one obtains that the fields $\Big(
X_a\Big)_{a=1,..., {\rm codim}C}$ are in {\em involution}:
$$
[X_a,X_b]=0.
$$
In other words, by Frobenius' theorem,  $(T C)^{\perp}$ is a
{\em integrable distribution}, the {\em characteristic} distribution. Notice that the
involution condition of the Hamiltonians
$$
\{H_a,H_b\}=0
$$
is intrinsic, constitutive, of the coisotropic submanifold $C$.
\medskip
\par\noindent
\section{Appendix: Geometric Theory of Ha\-milton-Jacobi equations on coisotropic submanifolds}
\medskip
\par\noindent
Let us consider $d$  functionally independent and in involution Hamiltonians ${\tilde
H}_a$,
$a=1,...,d$, of the following structure:
$$
{\tilde H}_a :T^*(N\times {\mathbb R}^d)\longrightarrow{\mathbb R}
$$
$$
(x_1,...x_n,t_1,,...,t_d,p_1,...p_n,\tau_1,,...,\tau_d)\longmapsto\tau_a +H_a(t_1,,...,t_d,x,p)
$$
The submanifold
$$
C=\bigcap_{a=1}^d{\tilde H}_a^{-1}(0)
$$
is coisotropic (because $\{ {\tilde H}_a , {\tilde H}_b \}=0$), ${\rm dim}\, C=2n+d$, ${\rm co}$-${\dim}\, C=d$.
\par\noindent
We call {\em geometric solution} of the Ha\-milton-Jacobi multi-equation related to
$C$ every {\em Lagrangian submanifold}\,\footnote{ that is, $i$)
$\sigma|_{\Lambda}=0$ and
$ii$)
${\rm dim}\,
\Lambda={\rm dim}\, N\times {\mathbb R}^d=n+d$.} $\Lambda\subset T^*(N\times {\mathbb R}^d) $
belonging to $C$:
$$
\Lambda\subset C.
$$
The  Lagrangian submanifolds which are graphs of a section of the cotangent bundle
$T^*(N\times {\mathbb R}^d)$ are image of the differential of functions $S$, defined on
$N\times {\mathbb R}^d$,\  \ \ $S:N\times {\mathbb R}^d\to {\mathbb R}$,
$$
\Lambda ={\rm im}\, dS = \Big\{(x,t,p,\tau): p=\frac{\partial S}{\partial x}(x,t),\  \tau=\frac{\partial S}{\partial t}(x,t)
\Big\}
$$
Thus, every globally transverse Lagrangian submanifold $\Lambda$, solving the H-J multi-equation related to $C$, represents, by its
generating function $S$,  also a  {\em classic} solution of it:
\begin{equation}\tag{A1}
\frac{\partial S}{\partial t_a}(x,t) +H_a(t_1,,...,t_q,x,\frac{\partial S}{\partial x}(x,t))=0,\quad a=1,..., d.
\end{equation}
\par\noindent
$\bullet$\ The characteristics method.
 The following proposition suggests us how to build geometric solutions  of the H-J multi-equation:
\par\noindent
{\bf Proposition}\ \ {\em Let $\Lambda\subset C$, $\Lambda$ be Lagrangian and $C=\cap_{a=1}^d{\tilde H}_a^{-1}(0)$
coisotropic. Then the Hamiltonian  vector fields  $X_{{\tilde H}_a}$, $a=1,..., d$, are tangent to $\Lambda$.
}
\par\noindent
 \begin{proof}If $v\in T_x\Lambda$ $\Rightarrow$  $v\in T_x C$ $\Rightarrow$ $v\in \cap_{a=1}^d{\rm ker}\, d\, {\tilde H}_a (x)$,
that is $d\, {\tilde H}_a (x)v=0,\ \ \forall a=1,...,d$. It holds that
$\sigma (X_{{\tilde H}_a},v)=-d\, {\tilde H}_a v=0 $, that is, the fields  $X_{{\tilde H}_a}$ are skew-orthogonal to $\Lambda$,
$X_{{\tilde H}_a}(x)\in (T_x\Lambda )^{\perp}= T_x\Lambda $, where the last equality holds since
$\Lambda
$
is Lagrangian.
 \end{proof}
\par\noindent
$\bullet$\ Multi-valued Variational Solutions.  The Liouville 1-form of $T^*(N\times  {\mathbb R}^d)$ is
$$\theta=\sum_{j=1}^n
p_j d x^j+\sum_{a=1}^d \tau_a d t^a$$ The restriction (pull-back) of $\theta$ on the Lagrangian submanifolds gives closed 1-forms,
 indeed, by denoting $j$
the inclusion map
$j:\Lambda \hookrightarrow T^*(N\times  {\mathbb R}^d)$,
$$
0=\sigma|_{\Lambda}=j^*\sigma=j^*d\theta=d j^*\theta=d(\theta|_{\Lambda}).
$$
Let $f:N\to {\mathbb R}$ be the initial data for the multi-equation (A1), then a primitive function of
$\theta|_{\Lambda}$ (on $\Lambda$, parametrized by $\chi, t^1,\dots ,t^d)$) is given by
$$
\bar S(\chi,t^1,...,t^d) =f(\chi)+\int_{
\tilde\Phi^{\bar t^1} \dots \tilde\Phi^{\bar t^a}  \dots  \tilde\Phi^{\bar t^d}(z)\big|_{\bar t^a\in [0,t^a]}}
\sum_{j=1}^n p_jdx^j-\sum_{a=1}^d H_a d t^a,
$$
$$
z=\left(\chi,0;df\left(\chi\right),-H_a\left(0,\chi,df\left(\chi\right)\right)\right)\in T^*(N\times {\mathbb R}^d)
$$
For small $t^1, \dots, t^d$, the $x$-components of $\tilde\Phi^{ t^1} \dots \tilde\Phi^{ t^a}  \dots
\tilde\Phi^{ t^d}(z)$, i.e.
$$
x(\chi, t)={\rm pr}_N \pi_{N\times {\mathbb R}^d}\, \tilde\Phi^{ t^1} \dots \tilde\Phi^{ t^a}  \dots
\tilde\Phi^{ t^d}(z),
$$
admits inverse: $\chi=\chi (x,t)$. In such a (scarcely meaningful) case, the function
$$
S(x,t):=\bar S(\chi (x,t),t)
$$
represents the classic solution to the Cauchy problem.
\medskip
\par\noindent
\section{Appendix: A short account on the developments of the theory of the multi-\-equations of
Ha\-milton-Jacobi}
\medskip
\par\noindent
After a pioneering paper by Tulczyjew \cite{T}, the study Hamilton-Jacobi problems in general co-isotropic
submanifolds of co-dimen\-sion grea\-ter than one inside $T^*Q$, that is, extending the standard case
$\frac{\partial S}{\partial t}+H(t,x,\frac{\partial S}{\partial x})=0$, was further developed by  Benenti and  Tulczyjew in
\cite{Ben2} and then in \cite{Ben}, where, among other things, in this last paper first it was considered generating
functions with finite auxiliary parameters for Lagrangian submanifolds of
$T^*Q\times T^*Q$ representing symplectic relations (canonical transformations). Examples of this new framework are the integrable
Dirac systems reconsidered by \cite{Lichnerowicz}  and \cite{M-T}.
More recently, new examples of H-J multi-equations arose from economics \cite{Rochet}, one can see the bibliography in
\cite{Barles-Tourin}. The use of symplectic topology and the notion of variational solution for the H-J equation have been
introduced by Chaperon and Sikorav, and then developed in \cite{Viterbo-X}, where,
in particular, the uniqueness of such a solution is deduced by the essential global uniqueness  of the generating functions proved in
\cite{Viterbo-STAGGF}. We recall that these variational solutions coincide with the viscosity solutions in the case of $p$-convex
Hamiltonians \cite{Zhukovskaia}, but it is not true in general case, lacking a (expected) Markov condition (cf. \cite{Viterbo-X}, \cite{Viterbo-Ottolenghi}).


\begin{thebibliography}{AKEEEEKDDDDEEEEEEE}


\bibitem[Barles]{Barles}
G. Barles, \newblock{\em Solutions de viscosit\'e des equations de Hamilton-Jacobi.}
{ Mathématiques et applications, vol. 17,} Springer-Verlag, 1994.

\bibitem [Barles-Tourin]{Barles-Tourin} G. Barles, A. Tourin, \newblock{  Commutation properties of
semigroups for first order Hamilton-Jacobi  equations and application to multi-time equations}, \newblock {\em
Indiana  Univ. Math. Journal},  vol 50, (2001), pp. 1524-1543.


\bibitem [Benenti]{Ben} S.Benenti, {\em Symplectic relations in analytical mechanics}, Proceedings of the
IUTAM-ISIMM symposium on modern developments in analytical mechanics, Vol. I (Torino, 1982). Atti Accad. Sci.
Torino Cl. Sci. Fis. Mat. Natur. 117 (1983), suppl. 1, 39--91.

\bibitem [Benenti-Tulczyjew]{Ben2} S.Benenti, W.M.Tulczyjew, {\em The geometric meaning and globalization of
the Hamil\-ton-Jacobi method}, Differential geometric methods in mathematical physics (Proc. Conf.,
Aix-en-Provence/Salamanca, 1979), pp. 9--21,  Lecture Notes in Math., 836, Sprin\-ger, Berlin, 1980.

 \bibitem[Bernardi-Cardin]{Bernardi-Cardin}
O.~Bernardi, F.~Cardin,
\newblock{Minmax and viscosity solutions of Hamilton-Jacobi equations in the
convex case}
\newblock{\em Communications on Pure and Applied analysis}
Volume: 5, Number: 4  (2006), pp. 793-812.
 
 \bibitem[Borel]{Borel} A.~Borel,\newblock{\em Cohomologie des espaces localement compacts (d'après Leray).}
 \newblock  Springer Lecture Notes, vol. 2,  Springer-Verlag, Berlin-Heidelberg-New-York.

 \bibitem[Brunella]{Brunella}
M.~Brunella,
\newblock{On a theorem of Sikorav}
\newblock{\em Enseignement Mathématique}
Volume:37,(1991), pp. 83--87.


  \bibitem[Chaperon]{Chaperon} M. Chaperon,
 \newblock{Lois de conservation et géométrie symplectique.}
   \newblock{\em  C. R. Acad. Sci. Paris. Série I, Math.},  312 (1991), pp. 345-348.

\bibitem[Chekanov]{Chekanov} Yu. V. Chekanov,
\newblock{ Critical points of quasifunctions, and generating families of Legendrian manifolds.} \newblock{ \em  Funktsional. Anal. i Prilozhen.} vol. 30 (1996),  pages 56--69.
\newblock{ Translation in Funct. Anal. Appl.vol.  30, (1996), no. 2, pages 118--128.}

\bibitem[Crandall-Lions]{Crandall-Lions}
M.G. Crandall and P.-L. Lions. 
\newblock{Viscosity solutions of Hamilton--Jacobi equations. }
\newblock{\em Trans. Amer. Math. Soc.}, 277:1--43, 1983.

\bibitem[Dacorogna-Marcellini]{Dacorogna-Marcellini} B. Dacorogna, P. Marcellini,
\newblock {\em Implicit Partial Differential equations.}
\newblock Birkhauser, Boston, 1999.

\bibitem[Ekeland-Hofer]{Ekeland-Hofer} I. Ekeland , H. Hofer, 
\newblock {Symplectic topology and Hamiltonian dynamics.}
\newblock {\em Math. Z.} 200 (1990) 355-378.

\bibitem[Eliashberg]{Eliashberg} Y. Eliashberg,
\newblock {Rigidity for symplectic and contact structures.}
\newblock  Preprint, 1981.

\bibitem[Entov-Polterovich-Zapolsky]{Entov-Polterovich-Zapolsky}
M. Entov, L. Polterovich and F. Zapolsky,
\newblock {\em Quasi-morphisms and the Poisson bracket.}
\newblock math.SG/0605406

\bibitem[Fathi-Maderna]{Fathi-Maderna}
 A. Fathi, E. Maderna,
 \newblock { Weak KAM theorem on non compact manifolds. }
 \newblock {\em Nolinear
Differential Equations and Applications}, vol. 11, n. 3, p. 271-300,
2004.

\bibitem[Gromov]{Gromov}
Misha Gromov.
\newblock Pseudoholomorphic curves in symplectic manifolds.
\newblock {\em Invent. Math.}, 82(2):307--347, 1985.

\bibitem[Guillemin-Pollack]{Guillemin-Pollack}
V.~Guillemin, A.~ Pollack,
\newblock{\em Differential topology},
 Prentice-Hall Inc.,Englewood Cliffs, N.J., 1974.

\bibitem [Hofer]{Hofer} H.~Hofer,\newblock {On the topological properties of symplectic maps.}
 \newblock {\em Proceedings
of the Royal Society of Edinburgh}, vol.115 (1990), pp. 25-38. 


 \bibitem [Humilière]{Humiliere}
 V.~Humilière.
 \newblock{\em On some completions of the space of Hamiltonian maps.}
 \newblock math.SG/0511418 (submitted)
 
 \bibitem [Humilière 2]{Humiliere2}
 V.~Humilière.
 \newblock{\em Quasi representations of nilpotent Lie algebras are $C^0$ closed.}
 \newblock in preparation, 2006.

\bibitem [Laudenbach-Sikorav]{Laudenbach-Sikorav} F. Laudenbach and J.-C. Sikorav, \newblock { Persistance d'intersection avec la section nulle au cours d'une isotopie hamiltonienne dans un fibré cotangent}, {\em Inventiones Math.} vol. 82 (1985) pp. 349-357. 

\bibitem[Laudenbach-Sikorav-2]{Laudenbach-Sikorav-IMRN}
Fran{\c{c}}ois Laudenbach and Jean-Claude Sikorav.
\newblock Hamiltonian disjunction and limits of lagrangian submanifolds.
\newblock {\em Intern. Math. Res. Notes}, 4:161--168, 1994.


\bibitem [Lichnerowicz] {Lichnerowicz} A. Lichnerowicz,
Vari\'et\'e symplectique et dynamique associ\'ee \`a une sous-vari\'et\'e. (French)
{\em C. R. Acad. Sci., Paris, Sér. } A 280, 523-527, (1975).

\bibitem [Lions-Rochet]{Lions-Rochet}  P.-L. Lions and J.-C. Rochet,
 \newblock {\em  Hopf formulas and
multi-time Hamilton-Jacobi equations}, \newblock {\em Proc.  Amer. Math. Soc.} 96 (1986), pp. 79-84.


\bibitem[Motta-Rampazzo]{Motta-Rampazzo}
Monica~Motta and Franco~Rampazzo,
\newblock Nonsmooth multi-time hamilton-jacobi systems.
\newblock {\em Indiana Univ. Math. J.}, 55:1573--1614, 2006.

\bibitem [Menzio-Tulczyjew]{M-T}  M.R. Menzio and W.M.Tulczyjew,
\newblock { Infinitesimal symplectic relations and generalized Hamiltonian dynamics}.
\newblock  {\em Ann. Inst. Henri Poincaré, Nouv. Sér.}, Sect. A 28, 349-367 (1978).

\bibitem[Rampazzo-Sussmann]{Rampazzo-Sussmann} F. Rampazzo, H. Sussmann, 
\newblock {\em  Set-valued differentials and a nonsmooth version of Chow's theorem. }
 \newblock{Published in the Proceedings of the 40th IEEE Conference on Decision and Control; Orlando, Florida, December 4 to 7, 2001 (IEEE Publications, New York, 2001), Volume 3, pp. 2613-2618.}

\bibitem[Rochet]{Rochet}  J.-C. Rochet,  \newblock {The taxation principle and multi-time Hamilton-Jacobi Equation} \newblock {\em  J. of Math. Econ.}14 (1985), 113-128.

\bibitem[Sikorav 1]{Sikorav-pc} J.-C. Sikorav, \newblock {Personnal communication.} 
\newblock { (1990). }

\bibitem[Sikorav-2]{Sikorav} J.-C. Sikorav, \newblock {
Problèmes d'intersections et de points fixes en géométrie hamiltonienne.} 
\newblock {Comment. Math. Helv. 62, No.1, pp. 62-73 (1987). }

\bibitem[Sikorav-3]{Sik3}
Jean-Claude Sikorav.
\newblock Quelques propriétés des plongements lagrangiens.
\newblock {\em Mém. S.M.F. n°46, suppl. au Bull; S.M.F. 119 (1991).}, pages
  151--167, 1991.

 \bibitem[Théret]{Theret-these}  D.~Théret,  \newblock {\em
\newblock{\em Thèse de Doctorat}
\newblock Université de Paris 7, Denis Diderot, 1996}.

\bibitem[Theret 2]{Theret} D.~Théret, \newblock {\em
\newblock{ A complete proof of Viterbo's uniqueness theorem on generating functions.}
\newblock {\em Topology Appl.} 96, no. 3, 249-266  (1999)}.

\bibitem [Tulczyjew]{T} W.~M.~Tulczyjew, \newblock {\em Relations symplectiques et les
\'equations d'Hamilton-Jacobi relativistes}, \newblock C. R. Acad.
Sc. Paris 281 A, 545-547, 1975.

\bibitem[Viterbo 1]{Viterbo-STAGGF} C.~Viterbo, \newblock { Symplectic  topology as the geometry of generating
functions}, \newblock  {\em Mathematische Annalen}, 292, (1992), pp.
685--710.

\bibitem[Viterbo 2]{Viterbo-X} C.~Viterbo,\newblock {\em Solutions d'\'equations de Hamilton-Jacobi},\newblock  Séminaire X-EDP 1995, Palaiseau.

 \bibitem[Viterbo 3]{Viterbo-isop} C.~Viterbo, \newblock { Metric and isoperimetric problems
  in symplectic topology}, \newblock {\em J. Amer. Math. Soc.} 13  (2000),  no. 2,
  pp. 411--431.

 \bibitem[Viterbo 4]{Viterbo-IMRN} C.~Viterbo, \newblock {On the uniqueness of generating Hamiltonian for continuous limits of Hamiltonians flows}, \newblock {\em International Mathematics Research Notices 2006} Article ID 34028, 9 pages, 2006. doi:10.1155/IMRN/2006/34028.
 
\bibitem[Viterbo-Ottolenghi]{Viterbo-Ottolenghi} C. Viterbo, A.  Ottolenghi, \newblock {\em Solutions d'équations d'Hamilton-Jacobi et géométrie symplectique}, \newblock  Preprint, 1995 \hfill \break \newblock {http://math.polytechnique.fr/cmat/viterbo/Semi.HJ.dvi}.

\bibitem[Weinstein]{Weinstein-CBMS} A.~Weinstein,  \newblock {\em Lectures on symplectic geometry},
\newblock  CBMS Lecture Series, AMS, 1974.

\bibitem [Zhukovskaya]{Zhukovskaia} T. Zhukovskaya,
\newblock{\em Singularités  de minimax et solutions faibles d'{E}quations aux dérivées partielles.}
\newblock Thèse de Doctorat, Université de Paris 7, 1993.

\bibitem [Zhukovskaya 2]{Zhukovskaia-2} T. Zhukovskaya,
\newblock{\em Metamorphoses of the Chaperon-Sikorav weak solutions of Hamilton-Jacobi equations}
\newblock Thèse de Doctorat, Université de Paris 7, 1993.

\end{thebibliography}
 \end{document}